\documentclass[10pt,a4paper]{article}
\usepackage[utf8]{inputenc}
\usepackage[english]{babel}
\usepackage{amsmath}
\usepackage{amsfonts}
\usepackage{amssymb}
\usepackage{graphicx}
\usepackage{subfigure}
\usepackage{enumerate}
\usepackage{subfigure}
\usepackage{algorithm}
\usepackage{algorithmic}
\usepackage{amsthm}
\usepackage{bm}
\usepackage{pgfplots}
\usepackage[backend=biber]{biblatex}

\newcommand{\ud}[1]{\mathbf{#1}}

\theoremstyle{plain}
\newtheorem{lm}{Lemma}[section]
\newtheorem{ap}{Assumption}[section]
\newtheorem{thm}{Theorem}[section]
\newtheorem{prop}{Proposition}[section]
\newtheorem{cor}{Corollary}[section]
\newtheorem{rmk}{Remark}[section]

\theoremstyle{definition}
\newtheorem{defn}{Definition}[section]
\newtheorem{example}{Example}[section]

\title{Solving Stochastic Optimization by Newton-type methods with Dimension-Adaptive Sparse Grid Quadrature}

\author{Yuancheng Zhou}

\begin{document}

\maketitle

\begin{abstract}
%Stochastic optimisation problems minimise expectations of random cost functions. One fundamental question we need to answer before computation of such problem is how to choose between 'optimise then discretise' and 'discretise then optimise'. We prove that in fact the 'optimise then discretise' method is a kind of generalisation of the 'discretise then optimise' method for stochastic optimisation problems. Accurate quadrature methods are required to evaluate the objective, gradient or Hessian. Dimension-adaptive sparse grid quadrature shows high accuracy in computing integral with smooth integrand. It is a kind of generalisation of the classical sparse grid method, which refines different dimensions according to their importance. We show that the dimension-adaptive sparse grid quadrature has better performance in the optimise then discretise' method than the 'discretise then optimise' method. We focus on Newton-type methods as our optimisation solver. In this case, we show the condition when the convergence may breakdown and develop stopping criteria which is aimed to reduce extra computational cost.
Stochastic optimisation problems minimise expectations of random cost functions. We use 'optimise then discretise' method to solve stochastic optimisation. In our approach, accurate quadrature methods are required to calculate the objective, gradient or Hessian which are in fact integrals.  We apply the dimension-adaptive sparse grid quadrature to approximate these integrals when the problem is high dimensional. Dimension-adaptive sparse grid quadrature shows high accuracy and efficiency in computing an integral with a smooth integrand. It is a kind of generalisation of the classical sparse grid method, which refines different dimensions according to their importance. We show that the dimension-adaptive sparse grid quadrature has better performance in the optimise then discretise' method than the 'discretise then optimise' method.

\end{abstract}

\section{Introduction}
Stochastic optimisation is a useful tool in decision making and has many applications. The general form of an unconstrained stochastic optimisation problem is
\begin{equation}
\min_{u\in U}\ \operatorname{\mathbb{E}}[h(u,W)],
\end{equation}
where $W$ is a $d$ dimensional random vector which is defined on the probability space $(\Omega,\mathcal{B},\mathbb{P})$, $\mathcal{B}$ is the Borel $\sigma$-algebra and $\mathbb{P}$ is the associate probability measure. $U$ is a subset of $\mathbb{R}^n$ which contains all possible decisions.

%In this paper, we will assume that we have a known distribution. In this case, we can further write the problem as
%\begin{equation}
%\min_{u\in U} \int_{\Omega} h(u,\ud{\omega})p(W)\ d\ud{\omega}.
%\end{equation}
%where $p(\ud{\omega})$ is probability density function of the random vector $W$. In order to be more concise, we further denote 
%\begin{equation}\label{F}
%f(u,\ud{\omega})=h(u,\ud{\omega})p(W),\ and\ \ F(u)=\int_{\Omega} f(u,\ud{\omega})\ d\ud{\omega}.
%\end{equation}
%Thus the problem is reduced into following form
%\begin{equation}
%\min_{u\in U} F(u).
%\end{equation}

If the random vector $W$ subjects to a probability density $p(w)$ on $\mathbb{R}^d$ the objective is of the form
\begin{equation}\label{F}
F(u):=\operatorname{\mathbb{E}} \left[h(u,W)\right]=\int_{\mathbb{R}^d}{}
h(u,w)p(w)\, dw = \int_{\mathbb{R}^d}f(u,w)\, dw
\end{equation}
where $f(u,w) = h(u,w) p(w)$.

There are two categories of approaches to solve the stochastic optimisation problem. One is based on the idea of 'discretise then optimise' while the other is based on 'optimise then discretise'. When we apply the 'discretise then optimise' method(DTOM) to solve the stochastic optimisation problem, it turns out to be the so called scenario generation method. The main idea of this kind of methods is to first approximate the integrals $(\ref{F})$ by a quadrature rule $Q$,
\begin{equation}\label{sur}
F(u)\approx S(u)=Q(f(u,\cdot)).
\end{equation}  
then use $S(u)$ as a surrogate objective function and minimise it. The key of the scenario generation method is to find a good approximation to the original objective function. Monte Carlo(MC) and Quasi Monte Carlo(QMC) methods are successfully used in scenario generation for many applications.There are also many studies in convergence of the MC and QMC methods~\cite{dick2013high,rubinstein2016simulation,metropolis1949monte,niederreiter1992random,sloan1998quasi}. Recently, Michael, Sanjay and David developed the scenario generation via the sparse grid method (SGSG)~\cite{MSD1,MSD2}). The advantage of their approach, as shown in their paper, is the SGSG method converges faster than these scenario generation methods based on MC and QMC if the integrand in the objective function is smooth enough. They also showed the epi-convergence of the SGSG. However, there are also disadvantages of SGSG. When an integral is approximated using a sparse grid quadrature, weights of some grid points can be negative. The problem brought from these negative weights is that some important properties of the objective function, e.g. convexity, are no longer kept~\cite{zhou2018application}. The original convex objective function can be replaced by a non-convex surrogate function in SGSG method. This will bring more difficulties in computing the problem, e.g. from convex optimisation to non-convex optimisation, and analyse the performance of the whole algorithm. 

%Thus, they have to use some non-convex optimisation solvers to solve the surrogate problem ~\cite{MSD2}. A non-convex optimisation problem is usually much more difficult than a convex one. In general, it is an NP-hard problems and it is time consuming to find the minimiser of such problem.

Here we will study the alternative 'optimise then discretise' method (OTDM) and apply the idea to solve the stochastic optimisation problems. The basic idea of our approach is to solve the system of equations
\begin{equation}
\nabla F(u)=G(u)=0.
\end{equation}
There are many numerical methods can be used to solve the system of the non-linear equations. In order not to be too general, we focus on the Newton-type methods in this paper. For the integrals appear during the computation of the gradient $G$, we will apply more sophisticated dimension adaptive sparse grid method.

Newton's method and its variants have been widely used in solving nonlinear optimisation problems and nonlinear equations~\cite{K1,K2}. Hessian matrix needs to be computed in the original Newton's method while most its variants such BFGS and L-BFGS-B are gradient based method. The iterations of Newton-type method can be written into the following form
\begin{equation}\label{iter_1}
u_{p+1}=u_p-A_p^{-1}G(u_p),\ p=1,2,\dots
\end{equation}
where $A_p$ is an approximation to the Hessian Matrix. The convergence of the Newton-type methods are well studied ~\cite{K1}. Newton's method is quadratically convergent when the initial value is close enough to the minimiser. Other Newton-type methods(Quasi Newton methods) have lower convergence rate than Newton's method~\cite{K1}, however, they are used more frequently in practical computation since the Hessian matrix is not required to be computed and stored. These convergence theories can make sure the sequence $\left\{u_p\right\}$ generated by the $(\ref{iter_1})$ converges to a minimiser. However, in practice, one can only get a perturbed sequence $\left\{\bar{u}_p\right\}$ rather than the ideal sequence $u_p$. This is because of the existence of both rounding errors and truncated errors during the approximation to the function value $F$, gradient $G$ and the computation of iterations $(\ref{iter_1})$. It is shown in~\cite{YP,K1,DW} that if the error from an approximation of the gradient $G$ is sufficiently small, then the perturbed Newton-type method will produce a sequence $\left\{\bar{u}_p\right\}$ which converges to a minimiser $u^*$. This result also implies the 
Newton-type method can stop convergence if the error in the approximation of the gradient $G$ is large to some extent. However,  computing the gradient $G$ involves computing high dimensional integrals when $d$ is large. It is difficult to get a very accurate approximation in this case. Thus it is important for us to know when to stop the iterations in the Newton-type method. If we stop it too early, the solution can not achieve its best potential accuracy. If we stop too late, we will waste a great amount of computational cost. 
%Since each entries of $G$ is an integral in our problem, it is important for us to choose a proper resolution of the grid which used to approximate $G$. If very fine grid is used, the computational cost will be very expansive, especially for high dimensional problems. On the other hand, if we use very coarse grid, the breakdown of the convergence of Newton-type method will occur at the very beginning of the iteration $(\ref{iter_1})$ and we can not get an accurate approximation to the minimizer as a result.
We offer the stopping criterion for our method based on the error analysis mentioned in ~\cite{YP}. 

If we assume the integral and the derivative are interchangeable, then the gradient $G$ is
\begin{equation}
G(u)=\nabla F(u)=\int_{\Omega}\nabla f(u,w)\ dw.
\end{equation}
Each component of $G$ is an integral. As we have mentioned before, if these integrals are high dimensional, it will be very difficult to compute them even for moderate accuracy. This is so called 'curse of dimensionality'. However, for special function classes, such as functions which have bounded mixed derivatives, the sparse gird method ~\cite{MG2,MG3,JG} can mitigate the curse of dimensionality to a large extent. The performance of the sparse grid method can be further improved if we treat each dimension differently. Actually, in many applications, the importance of different dimensions are not equal. This property inspires the idea of the dimension-adaptive sparse grid \cite{M1,MG1,M2,BH,MH}. Unlike the classical sparse grid method which treats all the dimensions equally during the computation, the dimensional adaptive sparse grid method always refines the most 'important' dimension first. The performance is thus improved for those integrals with dimensions of different importance. Both the sparse grid method and the dimension-adaptive sparse grid method are used when we compute high dimensional integrals in the computation.

This paper is organised as follows. In section 2, we introduce basic concepts and results in optimisation which we need to use. In section 3, we review the sparse grid and dimension adaptive sparse grid method. The 'optimise then discretise' method and algorithms based on it are given in section 4. In section 5, we study the convergence of the sequence generated by the algorithm in section 4. In section 6, we develop the stopping criterion for our method. Finally, we show some high dimensional examples in the last section.
%The error bounds of these two methods are also discussed in this section. We study the perturbed sequence of Newton-type method and give some criteria on how to choose the size of the grid based on the error analysis.
%Finally, some examples and numerical results are provided in the last section.  

\section{Basic concepts of optimisation}
First we begin with introducing some basic concepts of optimisation. 
\begin{defn}
$u^{*}$ is a global minimizer of $F$ if $F(u)\ge F(u^{*})$ for all $u\in \mathbb{R}^n$.
\end{defn}

\begin{defn}
$F$ is convex when $F$ satisfies
$$
F(tu+(1-t)v)\le tF(u)+(1-t)F(v),\quad \forall u,v\in \mathbb{R}^n,\ t\in [0,1].
$$
$F$ is $\gamma$-strongly convex when there exists $\gamma>0$ such that
$$
F(tu+(1-t)v)\le tF(u)+(1-t)F(v)-\frac{1}{2}\gamma t(1-t)\Vert u-v\Vert_2^2,\quad \forall u,v\in \mathbb{R}^n,\ t\in [0,1].
$$
\end{defn}

We then make some smoothness assumptions on the cost function $F$ such that the global minimiser of the stochastic optimisation problem exists. 

\begin{defn}
A function $g:\mathbb{R}^p\rightarrow \mathbb{R}^q$ is Lipschitz continuous with constant $L>0$ if
$$
\Vert g(x)-g(y)\Vert_2\le L\Vert x-y\Vert_2,\quad \forall x,y\in \mathbb{R}^q.
$$
\end{defn}

\begin{ap}\label{assumption}
The function $F$ is continuously differentiable and $\nabla F$ is Lipschitz continuous with constant $L>0$. In this case, we call $F$ $L$-smooth.
\end{ap}
By using Taylor expansion, we have $\Vert \nabla^2 F(u)\Vert_2\le L$ if $F$ is $L$-smooth and twice continuously differentiable.

The following Lemma gives connection between convexity of a function and its smoothness.

\begin{lm}[~\cite{Nocedal2006numerical}]
If $F$ is continuously differentiable, then $F$ is convex if and only if $F$ lies on or above any tangent line:
$$
F(v)\ge F(u)+\nabla F(u)^T(v-u),\quad \forall u,v\in \mathbb{R}^n.
$$
Also, $F$ is $\gamma$-strongly convex if and only if
$$
F(v)\ge F(u)+\nabla F(u)^T(v-u)+\frac{\gamma}{2}\Vert v-u\Vert_2^2,\quad \forall u,v\in \mathbb{R}^n.
$$
If $F$ is twice continuously differentiable, then $F$ is convex if and only if $\nabla^2 F(w)$ is positive semidefinite for every $w\in \mathbb{R}^n$. Also, $F$ is $\gamma$-strongly convex if and only if $\nabla^2 F(w)\ge \gamma I$.
\end{lm}

By using this Lemma, we can show the existence and uniqueness of global minimizers for strongly convex functions.

\begin{thm}[~\cite{Nocedal2006numerical}]\label{unique_solution}
If $F:\mathbb{R}^n\rightarrow \mathbb{R}$ is continuously differentiable and strongly convex, then it has a unique global minimizer.
\end{thm}

If the function $f$ satisfies assumptions in Theorem $\ref{unique_solution}$, we can make sure the stochastic optimisation problem is well defined. Next, we consider the numerical solvers to solve the problem. In order not to be too general, we will use Newton-type methods as our solvers. According to the optimality condition, solving the stochastic optimisation problem is equivalent to solve
 the following system of equations 
%Theorem $\ref{th1},\ref{th2},\ref{th3}$ not only make sure there exists global minimizer, but also illustrate how to find the it. The idea is to find a solution $u^*$ of the non-linear/linear equation
\begin{equation}\label{main_eq}
G(u)=\nabla F(u)=0.
\end{equation} 

The Newton-type methods generate following sequence $\left\{u_p\right\}$ 
\begin{equation}\label{iter}
u_{p+1}=u_p-A_p^{-1}G(u_p),\quad p=0,1,2,\dots\\
\end{equation}
and one expects the limit of this sequence will be the solution of $(\ref{main_eq})$. In the iteration, $A_p\in \mathcal{L}(U)$ is an 
approximation to the derivative $G^{'}(u)$, namely, the Hessian of $F$.
$A_p$ can be generated in many different ways and different choices of $A_p$ lead to different kinds of Newton-type methods. For example, if we take 
$$
A_p=\nabla^2F(u_p),
$$
this is exactly Newton's method. If we take 
$$
A_p=\alpha_p^{-1}B_p,
$$
where $\alpha_p$ is chosen by exact/inexact line search and $B_p$ is updated from the previously computed value $B_{p-1}$
$$
B_p=B_{p-1}+\frac{y_{p-1}y_{p-1}^T}{s_{p-1}^Ty_{p-1}}-\frac{(B_{p-1}s_{p-1})(B_{p-1}s_{p-1})^T}{s_{p-1}^TB_{p-1}s_{p-1}},
$$
where $B_0=I$, $s_{p}:=u_{p+1}-u_{p}$ and $y_{p}:=\nabla F(u_{p+1})-\nabla F(u_{p})$, then the iteration $(\ref{iter})$ becomes the BFGS method~\cite{K1,K2}, one of the most frequently used Quasi Newton methods.

If we further assume $\nabla^2F(u)$ is definite, bounded and Lipschitz continuous, see details in ~\cite{K1,K2}, the Newton method is quadratically convergent when the initial value is close enough to the minimiser while the BFGS method is superlinearly convergent. We can similarly assume the positive definiteness and boundedness of $\nabla^2f(u,w)$, $\forall w\in \Omega$ and the Lipschitz continuity of $\nabla^2f(u,w)$, $\forall w\in \Omega$ to make sure Newton method and BFGS method are convergent.

When we consider high dimensional problems, the difficulty lies in the approximations of $F(u)$, $\nabla F(u)$ and $\nabla^2 F(u)$.The objective $F(u)$ and each component of $\nabla F(u)$ and $\nabla^2 F(u)$ are high dimensional integrals in such case. Moreover, we have to compute them at each iteration in the solvers. This will result in the curse of dimensionality.

\section{Dimension-adaptive sparse grid}
\subsection{Formulation of the Dimension-adaptive sparse grid}
Here we introduce the sparse grid and the dimension-adaptive sparse grid method to approximate the $d$ dimensional integral
\begin{equation}
\int_{[-1,1]^d}f(x)\, dx.
\end{equation}
The sparse grid quadrature is built upon 1D quadrature. Suppose we have a sequence of 1D quadrature rules
\begin{equation}
Q_i(f)=\sum_{x_{i,j}\in G_i}c_{i,j}f(x_{i,j})
\end{equation}
where $G_i$ is a set which contains all quadrature points $\left\{x_{i,j},\  j=1,\dots, N_i\right\}$ of the $j$th 1D quadrature rule. $N_i$ is the number of the quadrature points. $\left\{w_{i,j},\  j=1,\dots, N_i\right\}$ are the corresponding weights of the $j$th rule. In particular, we focus on hierarchical quadrature rules here. The hierarchical means the sets $G_i$ are nested, i.e. $G_i\subset G_{i+1},\ \forall i$.

Once we have a sequence of 1D quadrature rules, we can define the differences
\begin{equation}
\Delta_i(f)=Q_i(f)-Q_{i-1}(f),\quad i=1,2,\dots
\end{equation}
with $Q_0(f)=0$. By using these differences, the $l$th quadrature rule can be written as
\begin{equation}
Q_l(f)=\sum_{i=1}^l\Delta_i(f).
\end{equation}

Multidimensional quadrature rules can be constructed based on the similar idea. Product rule is one of the methods used to deal with multidimensional integrals. The d-dimensional product rule is of the form
\begin{equation}\label{product_rule}
Q_{\ud{l}}(f)=\sum_{x_{\ud{l},\ud{j}}\in G_{\ud{l}}}c_{\ud{l},\ud{j}}f(X_{\ud{l},\ud{j}})
\end{equation}
where
\begin{equation}
\begin{aligned}
G_{\ud{l}}&=G_{l_1}\times\dots\times G_{l_d}\\
c_{\ud{l},\ud{j}}&=c_{l_1,j_1}\dots c_{l_d,j_d}\\
x_{\ud{l},\ud{j}}&=(x_{l_1,j_1},\dots,x_{l_d,j_d}).
\end{aligned}
\end{equation}
If we defined the d-dimensional differences as the tensor product of the 1D differences
\begin{equation}
\Delta_{\ud{i}}(f)=\Delta_{i_1}\otimes\dots\otimes\Delta_{i_d}(f),
\end{equation}
then the product rule $(\ref{product_rule})$ can be written as
\begin{equation}
Q_{\ud{l}}(f)=\sum_{\ud{i}\le \ud{l}}\Delta_{\ud{i}}(f).
\end{equation}

Instead of using all the differences in the index set $\left\{\ud{j},\ud{j}\le \ud{l}\right\}$, the sparse grid quadrature only sums over a subset of it. The level $l$ sparse grid quadrature is then defined by
\begin{equation}
Q_l^d(f)=\sum_{\ud{i}\le l+d-1}\Delta_{\ud{i}}(f).
\end{equation}

In fact, we can build a multidimensional quadrature by summing up any downset $\ud{I}$ of the full grid index set $\left\{\ud{i},\ \ud{i}\le \ud{l}\right\}$. The downset is defined as below.
\begin{defn}
We say $\ud{I}$ is a downset if it satisfies
\begin{equation}
\ud{i}\in \ud{I}\quad and \quad \ud{m}\le \ud{i},
\end{equation}
then $\ud{m}\in \ud{I}$.
\end{defn}

We call the multidimensional quadrature defined on such downset as generalised sparse grid quadrature and 
\begin{equation}
Q_{\ud{I}}^d(f)=\sum_{\ud{i}\in \ud{I}}\Delta_{\ud{i}}(f). 
\end{equation}
For the generalised sparse grid quadrature, the downset $\ud{I}$ is chosen before we compute the quadrature. Different downsets $\ud{I}$ are chosen in different applications. The truncated sparse grid, the sparse grid with fault, etc are several frequently used generalised sparse grid.

The dimension-adaptive sparse grid is still with the form
\begin{equation}
Q_{\ud{I}}^d(f)=\sum_{\ud{i}\in \ud{I}}\Delta_{\ud{i}}(f).
\end{equation}
while the downset $\ud{I}$ is decided during the computation according to the 'importance' of each dimension. There are two important things needed to be considered before designing the dimension adaptive algorithm. First, when we add a new surplus $\Delta_{\ud{i}}f$ to the sum, we need to make sure the newly generated index set $\mathbf{I}\cup \left\{\ud{i}\right\}$ is still a downset. This is because we need to use the method of differences to compute the telescope sum and thus every index $\ud{j}$ which have smaller entries than $\ud{i}$ in at least one dimension must be included in I. Second, the algorithm is required to detect the 'important dimension' and do refinement first in 'the most important dimension' at each iteration.

\begin{algorithm}
\caption{Dimension adaptive sparse grid quadrature}\label{dasg_alg}
\begin{algorithmic}
\STATE Initialize $\mathbf{I}=\left\{\ud{1}\right\}$ and $s=\Delta_{\ud{i}}f$
\WHILE{Termination condition not reached}

\STATE Consider all possible covering elements to $\mathbf{I}$ and put them in a heap $\mathcal{A}$
\STATE Select $\ud{i}$ from heap $\mathcal{A}$ with largest $\Delta_{\ud{i}}f$
\STATE s=s+$\Delta_{\ud{i}}f$

\ENDWHILE
\end{algorithmic}
\end{algorithm}
%
%The algorithm $1$ shows the basic idea of implementing dimension adaptive sparse grid. Here, we say $\left\{\ud{i}\right\}$ is a covering element of a downset $\mathbf{I}$ if $\mathbf{I}\cup \left\{\ud{i}\right\}$ is also a downset. By the definition of downset, judge whether $\left\{\ud{i}\right\}$  is a covering element is equivalent to check $\forall\ud{j}\le \ud{i}$, if $\ud{j}\in \mathbf{I}$ or not. Actually we can only check if all $\ud{j}=\ud{i}-e_j\in \mathbf{I}$ where $e_j$ is the $j$-th unit vector. This is because all $\ud{k}\le\ud{j}$ must be in the index set $\mathbf{I}$ since $\mathbf{I}$ is a downset. An efficient way to do this kind of checking is mentioned in $(\cite{MG1})$. The brief idea is to store a forward neighbors list and a backward neighbors list using the relative addresses of all index $\ud{i}$. When the first time the algorithm computes the surplus $\Delta_{\ud{i}}f$, a function is written to fill these two neighbors list for index $\ud{i}$. By comparing the address of $\ud{j}=\ud{i}-e_j\in \mathbf{I}$ with all the elements in the corresponding neighbors list of $\ud{i}$, we can judge whether $\ud{i}$ is a covering element of $\mathbf{I}$.
%
%We build a reference type heap to store all the addresses of the corresponding covering element $\ud{i}$ at current stage. The reference is also stored as a list in the heap. The value of the reference is the surplus $\Delta_{\ud{i}}$. The reference type heap is efficient to push the relative address of $\ud{i}$ in and pop the address of $\ud{k}$ associated with largest $\Delta_{\ud{k}}f$.

The termination condition we considered in this paper is
\begin{equation}\label{term}
\mathbf{I}=\left\{\ud{i}\le \ud{l}\ |\ \vert \Delta_{\ud{i}}f\vert\ge \epsilon\right\}.
\end{equation}
The termination condition $\vert \Delta_{\ud{i}}f\vert\geq \epsilon$ has been used in many dimension adaptive sparse grid algorithms to stop the while loop \cite{M2,BH,MH}. The additional condition $\ud{i}\le \ud{l}$ we added here is aimed at avoiding excessive refinement in some dimensions. Here, we say $\left\{\ud{i}\right\}$ is a covering element of a downset $\mathbf{I}$ if $\mathbf{I}\cup \left\{\ud{i}\right\}$ is also a downset.

We use $\mathcal{Q}_{\ud{l},\epsilon}$ as the operator for the dimension-adaptive sparse grid quadrature in Algorithm $\ref{dasg_alg}$ with the termination condition $(\ref{term})$. The choice of the downset  $\ud{I}$ depends on $f$, $\ud{l}$, $\epsilon$, so we have
\begin{equation}
\mathcal{Q}_{\ud{l},\epsilon}(f)=Q_{\ud{I}(f,\ud{l},\epsilon)}(f).
\end{equation}
It should be noted here we have to use different notations for the quadrature method($\mathcal{Q}_{\ud{l},\epsilon}$) and the computing formula($Q_{\ud{I}}$). This is because when the same quadrature method applied to approximate different integrals, e.g. integrals with integrand $f$ and $g$, respectively, we can get different downsets $\ud{I}_f$ and $\ud{I}_g$,
\begin{equation}
\ud{I}_f=\ud{I}(f,\ud{l},\epsilon)\neq\ud{I}(g,\ud{l},\epsilon)=\ud{I}_g.
\end{equation}
Thus, the formula used to approximate the integrals are different, that is,
\begin{equation}
Q_{\ud{I}_f}\neq Q_{\ud{I}_g}.
\end{equation}
For non-adaptive approach, we don't have such problem. We use the same notation($Q$) to denote both quadrature method and computing formula.

\subsection{1D quadrature rules}
We use the following three types of 1D quadrature rules to build our dimension-adaptive sparse grid quadratures. They are the trapezoidal rule, the Clenshaw-Curtis rule~\cite{CC} and the Gauss-Patterson rule~\cite{P1}. All of these are hierarchical quadrature rules. The trapezoidal rule has $O(4^{-l})$ accuracy on the uniform grid with $N=2^{l-1}+1$ grid points when the integrand $f\in C^2$. The accuracy can be further improved to $O(2^{-lr})$ if the integrand $f\in C^r$ is periodic. The $N=2^{l-1}+1$ points Clenshaw-Curtis rule uses extremal points of the Chebyshev polynomial of degree $N+1$ as its quadrature points. The $N+1$ points Clenshaw-Curtis rule integrates polynomials of degree less or equal than $N$ exactly~\cite{MG2,MH}. The accuracy of an $N$ points rule is $O(2^{-lr})$~\cite{DP,MG2}. The Gauss-Patterson rule is a Kronrod extension of the corresponding Gauss rule. The polynomial degree of exactness of an $N=2^{l}-1$ points rule is $(3N-1)/2$. Its accuracy for the integrals with integrand $f\in C^r$ is also $O(2^{-lr})$ ~\cite{DP,MG2} for an $N$ points rule. It is noteworthy that both level $l$ trapezoidal rule and Clenshaw-Curtis rule have $2^{l-1}+1$ quadrature points while the level $l$ Gauss-Patterson rule has $2^{l}-1$  quadrature points.

\subsection{Accuracy of high dimensional quadrature rules}
We will mainly discuss the accuracy of the dimension-adaptive sparse grid quadrature rules. Before that, we first provide the results on product rule and sparse grid quadrature. 

The computational complexity of the product rule is $O(N_l^d)$ for $l_i=l$. However, the accuracy is $O(2^{-lr})$. Here we notice that the accuracy is not depend on the dimension which results in the curse of dimensionality.

For the sparse grid quadrature rules, If we assume the integrand $f$ has bounded mixed derivatives up to order $r$, i.e., $f\in H^r([-1,1]^d)$ where
\begin{equation}
H^r([-1,1]^d)
=\left\{f:[-1,1]^d\rightarrow\mathbb{R}:\max_{\vert\bm{\alpha}\vert_{\infty}\le r}\left\|\frac{\partial^{\vert\bm{\alpha}\vert_1}f}{\partial^{\bm{\alpha}}W}\right\|<\infty\right\}
\end{equation}
where $\Vert\cdot\Vert$ denotes the $L_2$ norm and $\vert \bm{\alpha}\vert_{\infty}=\max_{j}\alpha_j$, then the error of the sparse grid quadrature is $O(N^{-r}(\log N)^{(d-1)(r-1)})$~\cite{MG2,MG3,MH}. The $N$ here is the number of the sparse grid quadrature points.

Next, we study the accuracy of the dimension-adaptive sparse grid quadrature. The error of the dimension-adaptive sparse grid interpolation has been studied in~\cite{M1,M2,BH}. Similar as the analysis in~\cite{M1}, we first have following bound on $\lvert Q_{L}f-Q_{I}f\rvert$.

%
%Next, we will look into the error bounds of $|Q_{\ud{l}}f-Q_{\mathbf{I}}|$. First, similar as the error bound in functional space developed in $(\cite{M1})$, we have following priori error bound for dimension adaptive sparse grid quadrature under termination condition $(\ref{term})$.

\begin{prop}\label{priori err}(Priori error bound)
Let $\mathbf{I}=\left\{\ud{i}\ |\ \ud{i}\le \ud{l},\ \vert \Delta_{\ud{i}}f\vert\ge \epsilon\right\}$ and $Q_{\ud{L}}f-Q_\mathbf{I}f$ be the error of the dimension adaptive sparse grid quadrature on set $\mathbf{I}$ relative to $Q_{\ud{l}}f$. We further denote $\mathbf{L}=\left\{\ud{i}\ |\ \ud{i}\le\ud{l}\right\}$ and the number of all indices in $\mathbf{L}$ is $\vert\mathbf{L}\vert=\prod_{k=1}^d(l_k+1)$. Then we get the bound
\begin{equation}
\vert Q_{\ud{L}}f-Q_\mathbf{I}f\vert\le \vert\mathbf{L}\vert\epsilon.
\end{equation}
\end{prop}
\begin{proof}
According to the definition, we have
\begin{equation}
\begin{aligned}
\vert Q_{\ud{L}}f-Q_\mathbf{I}f\vert &=\vert\sum_{\ud{i}\in \mathbf{L}}\Delta_{\ud{i}}f-\sum_{\ud{i}\in \mathbf{I}}\Delta_{\ud{i}}f\vert
=\vert \sum_{\ud{i}\in\mathbf{L}\backslash\mathbf{I}}\Delta_{\ud{i}}f\vert\\
&\le  \sum_{\ud{i}\in\mathbf{L}\backslash\mathbf{I}}\vert\Delta_{\ud{i}}f\vert\le\sum_{\ud{i}\in\mathbf{L}\backslash\mathbf{I}}\epsilon\le\sum_{\ud{i}\in\mathbf{L}}\epsilon= \vert\mathbf{L}\vert\epsilon.
\end{aligned}
\end{equation}
The first inequality follows from the triangle inequality. The second inequality holds because $\mathbf{L}\backslash\mathbf{I}=\left\{\ud{i}\le \ud{l}\ |\ \vert \Delta_{\ud{i}}f\vert< \epsilon\right\}$. The third inequality follows by the fact $\mathbf{I}\subset\mathbf{L}$.
\end{proof}

From the proof of the proposition $(\ref{priori err})$, we do not use any information from the computation process of the $Q_{\mathbf{I}}f$. The error bound can be derived before computing $Q_{\mathbf{I}}f$. However, after we compute $Q_{\mathbf{I}}f$ by the dimension-adaptive sparse grid method, we will know exactly what the downset $\mathbf{I}$ is. This can help us improve the error bound.

\begin{prop}\label{posteriori bound}(Posteriori error bound)
Suppose the downset $\mathbf{I}$ is known after we computed the $Q_{\mathbf{I}}f$. The error bound in $\ref{priori err}$ can be improved by
\begin{equation}
\vert Q_{\ud{L}}f-Q_\mathbf{I}f\vert\le (\vert\mathbf{L}\vert-\vert\mathbf{I}\vert)\epsilon,
\end{equation}
where the set $\mathbf{L}=\left\{\ud{i}\ |\ \ud{i}\le\ud{l}\right\}$.
\end{prop}
\begin{proof}
Since $\mathbf{I}=\left\{\ud{i}\ |\ \ud{i}\le \ud{l},\ \vert \Delta_{\ud{i}}f\vert\ge \epsilon\right\}$ is a subset of $\mathbf{L}$, we have
\begin{equation}
\vert Q_{\ud{L}}f-Q_\mathbf{I}f\vert=\vert\sum_{\ud{i}\in \mathbf{L}\backslash\mathbf{I}} \Delta_{\ud{i}}f\vert\le \sum_{\ud{i}\in \mathbf{L}\backslash\mathbf{I}} \vert\Delta_{\ud{i}}f\vert=(\vert \mathbf{L}\vert-\vert \mathbf{I}\vert)\epsilon. 
\end{equation}

\end{proof}

When we compute a high dimensional integral, we won't set a very small  $\epsilon$, e.g. $10^{-15}$, in the termination condition since the computational cost is usually unaffordable for most cases. Thus, neither priori error bound nor posteriori error bound is good when we consider a high dimensional problem since $\vert\mathbf{L}\vert$ will grow exponentially when the dimension increases while the $\epsilon$ can't be chosen as small as possible. In order to get a more accurate error bound, we need to use the smoothness of the integrand $f$. 

\begin{lm}\label{bound for d difference operator}
If $f\in H^r([-1,1]^d)$ and we have following estimation for the 1D quadrature rules $Q_i$ 
\begin{equation}
\vert If-Q_if\vert\le \Vert I-Q_i\Vert\Vert f\Vert\le \gamma_r2^{-ri}\Vert f\Vert
\end{equation}
where the constants $\gamma_r$ can be obtained by known bounds for the respective Peano kernels(ref),
then if hierarchical rules are used in building $d$ dimensional quadrature rules, we have $\vert\Delta_{\ud{i}}f\vert\le C_r^d2^{-r\vert\ud{i}\vert}\Vert f\Vert$ where $C_r^d=\gamma_r^d(1+2^r)^d$.
\end{lm}
%\begin{ap}\label{surplus assump}
%We assume the surplus with index $\ud{i}$ satisfies $\vert\Delta_{\ud{i}}f\vert\le K 4^{-\vert\ud{i}\vert}$ where $K$ is some constant.
%\end{ap}
\begin{proof}
Since we use hierarchical rules, the 1D difference can be written as
\begin{equation}
\begin{aligned}
\Delta_i f=Q_i f-Q_{i-1}f&=\sum_{x_{i,j}\in G_i}c_{i,j}f(x_{i,j})-\sum_{x_{i-1,j}\in G_{i-1}}c_{i-1,j}f(x_{i-1,j})\\
&=\sum_{x_{i,j}\in G_i}b_{i,j}f(x_{i,j})
\end{aligned}
\end{equation}
where $b_{i,j}=c_{i,j}$ for quadrature points in the set $G_i\backslash G_{i-1}$ and $b_{i,j}=c_{i,j}-c_{i-1,j}$ otherwise. Then for 2D case, we have
\begin{equation}
\Delta_{i_1}\otimes\Delta_{i_2}f=\sum_{x_{i_1,j_1}\in G_{i_1}}\sum_{x_{i_2,j_2}\in G_{i_2}}b_{i_1,j_1}b_{i_2,j_2}f(x_{i_1,j_1},x_{i_2,j_2}),
\end{equation} 
Furthermore,
\begin{equation}
\begin{aligned}
\vert\Delta_{i_1}\otimes\Delta_{i_2}f\vert&\le \Vert\Delta_{i_1}\Vert\Vert\sum_{x_{i_2,j_2}\in G_{i_2}}b_{i_2,j_2}f(\cdot, x_{i_2,j_2})\Vert\\
&\le \Vert\Delta_{i_1}\Vert\sup_{0\le\alpha_1\le r}\sup_{s\in [-1,1]}\vert\Delta_{i_2}f^{(\alpha_1,0)}(s,\cdot)\vert\\
&\le \Vert\Delta_{i_1}\Vert\sup_{0\le\alpha_1\le r}\sup_{s\in [-1,1]} \Vert\Delta_{i_2}\Vert\Vert f^{(\alpha_1,0)}(s,\cdot)\Vert\\
&\le \Vert\Delta_{i_1}\Vert\Vert\Delta_{i_2}\Vert\Vert f\Vert.
\end{aligned}
\end{equation}
This can be generated to $d$ dimensional case that is
\begin{equation}\label{d tensor product estimation}
\vert\Delta_{\ud{i}}f\vert=\vert\Delta_{i_1}\otimes\Delta_{i_2}\otimes\dots\otimes\Delta_{i_d}f\vert\le \Vert\Delta_{i_1}\Vert \Vert\Delta_{i_2}\Vert\dots\Vert\Delta_{i_d}\Vert\Vert f\Vert.
\end{equation}
We can derive the upper bound of the norm of the 1d difference operator from 
\begin{equation}\label{1d tensor product estimation}
\begin{aligned}
\Vert\Delta_{i_k}\Vert&=\Vert Q_{i_k}-Q_{i_{k-1}}\Vert\\
&\leq\Vert I-Q_{i_k}\Vert+\Vert I-Q_{i_{k-1}}\Vert\\
&\leq r_{\gamma} 2^{-ri_k}(1+2^r).
\end{aligned}
\end{equation}
Combining $(\ref{d tensor product estimation})$ and $(\ref{1d tensor product estimation})$, we get
\begin{equation}
\vert\Delta_{\ud{i}}f\vert\le C_r^d2^{-r\vert\ud{i}\vert}\Vert f\Vert.
\end{equation}
\end{proof}

\begin{thm}\label{surplus assump bound}
Under the conditions of Lemma $(\ref{bound for d difference operator})$, we can further improve our posteriori bound by
\begin{equation}\label{surplus assump bound formula}
\vert Q_{\ud{L}}f-Q_\mathbf{I}f\vert\le K\sum_{\ud{i}\in \mathbf{L}\backslash\mathbf{I}} 2^{-r\vert\ud{i}\vert}.
\end{equation}

\end{thm}
\begin{proof}
The proposition $(\ref{surplus assump bound})$ follows directly from the proposition $(\ref{posteriori bound})$ and the Lemma $(\ref{bound for d difference operator})$. $K=C_r^d\lVert f\rVert$ is a constant.
\end{proof}

\begin{cor}\label{fin_err_bound}
Suppose $\ud{m}$ is one of the indices such that $\vert\ud{m}\vert\le \vert\ud{i}\vert,\forall \ud{i}\in \mathbf{L}\backslash\mathbf{I}$, then the posterior bound is
\begin{equation}\label{err_bound_ineq}
\vert Q_{\ud{L}}f-Q_\mathbf{I}f\vert\le \frac{\epsilon}{\rho}\sum_{\ud{i}\in \mathbf{L}\backslash\mathbf{I}}2^{r(\vert\ud{m}\vert-\vert\ud{i}\vert)}
\end{equation}
where $\rho K2^{-r\vert\ud{m}\vert}=\epsilon$ and 
\begin{equation}
\rho_{\min}:=\frac{\sum_{\ud{i}\in \mathbf{L}\backslash\mathbf{I}}2^{r(\vert\ud{m}\vert-\vert\ud{i}\vert})}{\vert \mathbf{L}\vert-\vert\mathbf{I}\vert}\le\rho< 2^r.
\end{equation}
\end{cor}
\begin{proof}
We first notice we can rewrite $(\ref{surplus assump bound formula})$ as
\begin{equation}
\vert Q_{\ud{L}}f-Q_\mathbf{I}f\vert\le K2^{-r\vert\ud{m}\vert}\sum_{\ud{i}\in \mathbf{L}\backslash\mathbf{I}}2^{r(\vert\ud{m}\vert-\vert\ud{i}\vert)}.
\end{equation}
By using $\rho K2^{-r\vert\ud{m}\vert}=\epsilon$, we get the inequality 
$(\ref{err_bound_ineq})$. For the lower bound of $\rho$, we expect the error bound $(\ref{err_bound_ineq})$ is not worse than the posterior error bound in the proposition $(\ref{posteriori bound})$,otherwise we can use the latter one. Thus, we have
\begin{equation}
\frac{\epsilon}{\rho}\sum_{\ud{i}\in \mathbf{L}\backslash\mathbf{I}}2^{r(\vert\ud{m}\vert-\vert\ud{i}\vert)}\le \epsilon(\vert \mathbf{L}\vert-\vert\mathbf{I}\vert),
\end{equation}
This leads to the lower bound of $\rho$. For the upper bound, if we denote the $k$th negative unit vector as $\ud{e}_k=[0,\dots,1,\dots,0]$, then according to the definition of $\ud{m}$, there exists an index $\ud{m}-\ud{e}_k\in \mathbf{I}$, otherwise we should choose $\ud{m}-\ud{e}_k$ instead of $\ud{m}$ in the theorem. According to the Theorem $\ref{surplus assump bound}$ and the definition of the downset $\mathbf{I}$, we have
\begin{equation}
\epsilon<\vert\Delta_{\ud{m}-\ud{e}_k}f\vert\le K 2^{-r(\vert\ud{m}\vert-1)}.
\end{equation}
Thus, using the definition of $\rho$, we have
\begin{equation}
\rho=\frac{\epsilon}{K2^{-r\vert\ud{m}\vert}}\le \frac{K 2^{-r(\vert\ud{m}\vert-1)}}{K2^{-r\vert\ud{m}\vert}}=2^r.
\end{equation}
\end{proof}

By using the error estimation in the Theorem $\ref{fin_err_bound}$ and the error bound for d-dimensional product rule, we can obtain a bound on $|If-Q_\mathbf{I}f|$ by using trapezoidal rule, that is
\begin{equation}\label{bound}
\vert If-Q_\mathbf{I}f\vert\le \vert If-Q_{\ud{L}}f\vert+\vert Q_{\ud{L}}f-Q_\mathbf{I}f\vert\le c_d2^{-lr}+K\sum_{\ud{i}\in \mathbf{L}\backslash\mathbf{I}} 2^{-r\vert\ud{i}\vert}
\end{equation}
%\frac{\epsilon}{\rho}\sum_{\ud{i}\in \mathbf{L}\backslash\mathbf{I}}4^{\vert\ud{j}\vert-\vert\ud{i}\vert}. 
In the above bound, isotropic grid $\mathcal{G}_{\ud{l}}$, $l_i=l$ is used in the comparison. In ~\cite{M2}, the authors get an optimised priori error bound for $|If-Q_\mathbf{I}f|$ by balancing error bounds for the term $\vert If-Q_{\ud{L}}f\vert$ and $\vert Q_{\ud{L}}f-Q_\mathbf{I}f\vert$. The $\epsilon$ need to be chosen very small in order to achieve the optimized bound for high dimensional problems. In this paper, we are more interested in the case when
\begin{equation}
\vert If-Q_{\ud{L}}f\vert\ll\vert Q_{\ud{L}}f-Q_\mathbf{I}f\vert
\end{equation}
which means the approximation of the integral need to be accurate to some extent on the corresponding full grid, otherwise we can not expect the dimensional adaptive sparse grid method which uses a subset of quadrature points on the full grid provides a good approximation. Larger $\epsilon$ is allowed in this situation.
%This bound gives some hints on how to choose a suitable $\epsilon$ in our termination condition based on the smoothness of the integrand, dimension of the integrand and the cost of computation. We can minimize the right-hand side of $(\ref{bound})$ with respect to $n$ for some given $\epsilon$. This will give us 
%\begin{equation}
%\epsilon_n=c_d(n+1)^{-d}2^{-nr},
%\end{equation}  
%and the corresponding error bound is $2c_d2^{-nr}$.(To be continued here, illustrate by plots).
\section{The 'optimise then discretise' method}
We will show the framework of the 'optimise then discretise' method and test a 2D example to illustrate its performance. For simplicity, we first use Newton method as our optimisation algorithm. It is shown in the Algorithm $\ref{optimise}$. The Algorithm $\ref{discretise}$, $\ref{general_discretise}$, $\ref{discretise(adaptive)}$, $\ref{modified_discretise(adaptive)}$ and $\ref{general_discretise(adaptive)}$ are discretised versions of the Newton method in Algorithm $\ref{optimise}$ from simple to complex. In the Algorithm $\ref{discretise}$ and $\ref{general_discretise}$, we use non-adaptive quadrature $Q$ to compute the integrals. In the Algorithm $\ref{modified_discretise(adaptive)}$ and $\ref{general_discretise(adaptive)}$, we use the dimension adaptive quadrature $\mathcal{Q}$ to compute the integrals. The notation $D_i$ denotes the $i$th discretised derivative and thus $D_{ij}^2$ is the $ij$th second order discretised derivative. $D$ denotes the discretised gradient and $D^2$ denotes the discretised Hessian.

In Algorithm $\ref{discretise}$, we use different non-adaptive surrogate at each iteration. The non-adaptive quadrature operator $Q$ and the discretised derivative operators $D_i$, $D_{ij}$ are commutative, i.e.
$$
D_iQ_p=Q_pD_i\ \text{and}\ D_{ij}^2Q_p=Q_pD_{ij}^2.
$$ 
This is because both two operators are fixed finite summations. In Algorithm $\ref{general_discretise}$, we further allow different choices of the non-adaptive quadrature for the objective and different component of gradient and Hessian.

The Algorithm $\ref{discretise(adaptive)}$ looks almost the same as the Algorithm $\ref{discretise}$ except the quadrature method is dimension adaptive. However, they have essential differences. The dimension adaptive operator $\mathcal{Q}$ and the discretised derivative operators $D_i$, $D_{ij}$ are not commutative. This observation leads to a new Algorithm $\ref{modified_discretise(adaptive)}$. The reason why the dimension adaptive operator and the discretised derivative operators are not commutative is because the downsets used in the computation are not equal, i.e.
$$
\begin{aligned}
I_{(D_if,\epsilon_p,\ud{l}_p)}&\neq I_{(f,\epsilon_p,\ud{l}_p)}\\
I_{(D_{ij}^2f,\epsilon_p,\ud{l}_p)}&\neq I_{(f,\epsilon_p,\ud{l}_p)}.
\end{aligned}
$$
One can also generalise the Algorithm $\ref{modified_discretise(adaptive)}$ to the Algorithm $\ref{general_discretise(adaptive)}$ by allowing the usage of different parameters in dimension adaptvie quadrature for the objective and different component of gradient and Hessian. Though Algorithm $\ref{general_discretise(adaptive)}$ will be more flexible than Algorithm $\ref{modified_discretise(adaptive)}$, we mostly use Algorithm $\ref{modified_discretise(adaptive)}$ in practice because it is usually hard to get information used for choose different parameters and the algorithm $\ref{general_discretise(adaptive)}$ is too complex.

%Its discretised versions are the Algorithm $\ref{discretise}$ and the Algorithm $\ref{discretise(adaptive)}$. It should be noted here there are essential differences between the adaptive discretisation and the non-adaptive one. Therefore, we write them separately.

\begin{algorithm} 
\caption{$\mathbf{OPTIMISE}$}\label{optimise}
\begin{algorithmic}[1]
\STATE Take an initial $u_0\in \mathbb{R}^n$ and $p:=0$
\STATE Compute $G_0=\nabla F(u_0)$
\WHILE{$\lVert G_p\rVert>\epsilon$}
\STATE Compute the Hessian $H_p=\nabla^2F(u_p)$
\STATE Update 
$$
u_{p+1}=u_p-H_p^{-1}G_p
$$
\STATE Set $p:=p+1$
\STATE Compute $G_p=\nabla F(u_p)$
\ENDWHILE
\STATE Output $u_p$ and $F(u_p)$
\end{algorithmic}
\end{algorithm}

\begin{algorithm}
\caption{$\mathbf{DISCRETISED\ VERSION}$}\label{discretise}
\begin{algorithmic}[1]
\STATE Take an initial $\bar{u}_0\in \mathbb{R}^n$ and $p:=0$
\STATE Compute the approximation of the gradient $\bar{G}_0=DQ_{0}(f(\bar{u}_0,\cdot))$
\WHILE{$\lVert \bar{G}_p\rVert>\epsilon$}
\STATE Compute the approximation of the Hessian $\bar{H}_p=D^2Q_p(f(\bar{u}_p,\cdot))$
\STATE Update
\begin{equation}\label{OTDM_iter}
\bar{u}_{p+1}=\bar{u}_p-\bar{H}_p^{-1}\bar{G}_p
\end{equation}
\STATE Set $p:=p+1$
\STATE Compute the approximation of the gradient $\bar{G}_p=DQ_{p}(f(\bar{u}_p,\cdot))$
\ENDWHILE
\STATE Output $\bar{u}_p$ and $\bar{F}_p:=Q_{p}(f(\bar{u}_p,\cdot))$
\end{algorithmic}
\end{algorithm}

\begin{algorithm}
\caption{$\mathbf{GENERAL\ DISCRETISED\ VERSION}$}\label{general_discretise}
\begin{algorithmic}[1]
\STATE Take an initial $\bar{u}_0\in \mathbb{R}^n$ and $p:=0$
\STATE Compute the approximation of the gradient $\bar{G}_0=[Q_{0,i}^{G}(D_if(\bar{u}_0,\cdot))]_{n\times 1}$
\WHILE{$\lVert \bar{G}_p\rVert>\epsilon$}
\STATE Compute the approximation of the Hessian $\bar{H}_p=[Q_{p,i,j}^{H}(D_{ij}^2f(\bar{u}_p,\cdot))]_{n\times n}$
\STATE Update
\begin{equation}\label{OTDM_iter}
\bar{u}_{p+1}=\bar{u}_p-\bar{H}_p^{-1}\bar{G}_p
\end{equation}
\STATE Set $p:=p+1$
\STATE Compute the approximation of the gradient $\bar{G}_p=[Q_{p,i}^{G}(D_if(\bar{u}_p,\cdot))]_{n\times 1}$
\ENDWHILE
\STATE Output $\bar{u}_p$ and $\bar{F}_p:=Q_{p}^{O}(f(\bar{u}_p,\cdot))$
\end{algorithmic}
\end{algorithm}

\begin{algorithm}
\caption{$\mathbf{DISCRETISED\ VERSION(ADAPTIVE)}$}\label{discretise(adaptive)}
\begin{algorithmic}[1]
\STATE Take an initial $\bar{u}_0\in \mathbb{R}^n$ and $p:=0$
\STATE Compute the approximation of the gradient $\bar{G}_0=D\mathcal{Q}_{\epsilon_{0}, \ud{l}_{0}}(f(\bar{u}_0,\cdot))$
\WHILE{$\lVert \bar{G}_p\rVert>\epsilon$}
\STATE Compute the approximation of the Hessian $\bar{H}_p=D^2\mathcal{Q}_{\epsilon_p, \ud{l}_p}(f(\bar{u}_p,\cdot))$
\STATE Update
\begin{equation}\label{OTDM_iter}
\bar{u}_{p+1}=\bar{u}_p-\bar{H}_p^{-1}\bar{G}_p
\end{equation}
\STATE Set $p:=p+1$
\STATE Compute the approximation of the gradient $\bar{G}_p=D\mathcal{Q}_{\epsilon_{p}, \ud{l}_{p}}(f(\bar{u}_p,\cdot))$
\ENDWHILE
\STATE Output $\bar{u}_p$ and $\bar{F}_p:=\mathcal{Q}_{\epsilon_{p}, \ud{l}_{p}}(f(\bar{u}_p,\cdot))$
\end{algorithmic}
\end{algorithm}

\begin{algorithm}
\caption{$\mathbf{MODIFIED\ DISCRETISED\ VERSION(ADAPTIVE)}$}\label{modified_discretise(adaptive)}
\begin{algorithmic}[1]
\STATE Take an initial $\bar{u}_0\in \mathbb{R}^n$ and $p:=0$
\STATE Compute the approximation of the gradient $\bar{G}_0=[\mathcal{Q}_{\epsilon_{0}, \ud{l}_{0}}(D_{i}f(\bar{u}_0,\cdot))]_{n\times 1}$
\WHILE{$\lVert \bar{G}_p\rVert>\epsilon$}
\STATE Compute the approximation of the Hessian $\bar{H}_p=[\mathcal{Q}_{\epsilon_p, \ud{l}_p}(D_{ij}^2f(\bar{u}_p,\cdot))]_{n\times n}$
\STATE Update
\begin{equation}\label{OTDM_iter}
\bar{u}_{p+1}=\bar{u}_p-\bar{H}_p^{-1}\bar{G}_p
\end{equation}
\STATE Set $p:=p+1$
\STATE Compute the approximation of the gradient $\bar{G}_p=[\mathcal{Q}_{\epsilon_{p}, \ud{l}_{p}}(D_if(\bar{u}_p,\cdot))]_{n\times 1}$
\ENDWHILE
\STATE Output $\bar{u}_p$ and $\bar{F}_p:=\mathcal{Q}_{\epsilon_{p}, \ud{l}_{p}}(f(\bar{u}_p,\cdot))$
\end{algorithmic}
\end{algorithm}

\begin{algorithm}
\caption{$\mathbf{GENERAL\ DISCRETISED\ VERSION(ADAPTIVE)}$}\label{general_discretise(adaptive)}
\begin{algorithmic}[1]
\STATE Take an initial $\bar{u}_0\in \mathbb{R}^n$ and $p:=0$
\STATE Compute the approximation of the gradient 
$$
\bar{G}_0=[\mathcal{Q}_{\epsilon_{0,i}, \ud{l}_{0,i}}^G(D_{i}f(\bar{u}_0,\cdot))]_{n\times 1}
$$
\WHILE{$\lVert \bar{G}_p\rVert>\epsilon$}
\STATE Compute the approximation of the Hessian 
$$
\bar{H}_p=[\mathcal{Q}_{\epsilon_{p,i,j}, \ud{l}_{p,i,j}}^H(D_{ij}^2f(\bar{u}_p,\cdot))]_{n\times n}
$$
\STATE Update
\begin{equation}\label{OTDM_iter}
\bar{u}_{p+1}=\bar{u}_p-\bar{H}_p^{-1}\bar{G}_p
\end{equation}
\STATE Set $p:=p+1$
\STATE Compute the approximation of the gradient 
$$
\bar{G}_p=[\mathcal{Q}_{\epsilon_{p,i}, \ud{l}_{p,i}}^{G}(D_if(\bar{u}_p,\cdot))]_{n\times 1}
$$
\ENDWHILE
\STATE Output $\bar{u}_p$ and $\bar{F}_p:=\mathcal{Q}_{\epsilon_{p}, \ud{l}_{p}}^{O}(f(\bar{u}_p,\cdot))$
\end{algorithmic}
\end{algorithm}

For more complicated quasi Newton methods, we take BFGS method with the exact line search as an example. The optimise algorithm and its adaptive discretised version are shown in the Algorithm $\ref{bfgs_optimise}$ and Algorithm $\ref{bfgs_discretise}$ which are extensions of the Algorithm $\ref{optimise}$ and the Algorithm $\ref{modified_discretise(adaptive)}$, respectively. In practice, the exact line search $(\ref{exact_line_search})$ is replaced with inexact line search for efficiency. The commonly used inexact line search is strong Wolfe's rule. The sequence $ u_p$ generated by BFGS with Wolfe's rule is proved to converge to the exact minimizer $ u^*$ superlinearly~\cite{K2}. It should be noted that we need to compute the objective $F$ in the line search methods(include the strong Wolfe's rule used in practice)in each iteration while this is not required if we use Newton method. 

\begin{algorithm} 
\caption{$\mathbf{BFGS\ OPTIMISE}$}\label{bfgs_optimise}
\begin{algorithmic}[1]
\STATE Take an initial $ u_0\in \mathbb{R}^n$ , an initial positive definite matrix $H_0$ and $p:=0$
\WHILE{$\lVert G( u_p)\rVert>\epsilon$}
\STATE Compute the search direction $ v_p=-H_pG( u_p)$
\STATE Find the step length $\alpha_p$ by exact line search
\begin{equation}\label{exact_line_search}
\min_{\alpha_p} F( u_p+\alpha_p v_p).
\end{equation}
The underlying $A_p^{-1}$ here is $\alpha_pH_p$. 
\STATE Update 
$$
 u_{p+1}= u_p+\alpha_p v_p
$$
\STATE Define $s_p:= u_{p+1}- u_p$ and $y_p:=G( u_{p+1})-G( u_p)$
\STATE Update
$$
H_{p+1}=\left(I-\frac{s_py_p^T}{s_p^Ty_p}\right)H_p\left(I-\frac{y_ps_p^T}{s_p^Ty_p}\right)+\frac{s_ps_p^T}{s_p^Ty_p}
$$
\STATE p:=p+1
\ENDWHILE
\STATE Output $ u_p$ and $F( u_p)$
\end{algorithmic}
\end{algorithm}

\begin{algorithm}
\caption{$\mathbf{BFGS\ DISCRETISE}$}\label{bfgs_discretise}
\begin{algorithmic}[1]
\STATE Take an initial $\bar{ u}_0\in \mathbb{R}^n$, an initial positive definite matrix $\bar{H}_0$ and $p:=0$
\STATE Compute $\bar{G}_0=[\mathcal{Q}_{\epsilon_{0}, \ud{l}_{0}}(D_{i}f(\bar{u}_0,\cdot))]_{n\times 1}$
\WHILE{$\lVert \bar{G}_p\rVert>\epsilon$}
\STATE Compute the search direction $\bar{ v}_p=-\bar{H}_p\bar{G}_p$
\STATE Find the step length $\bar{\alpha}_p$ by exact line search
$$
\min_{\bar{\alpha}_p} \bar{F}_p(\bar{ u}_p+\bar{\alpha}_p\bar{ v}_p)
$$
where $\bar{F}_p:=\mathcal{Q}_{\epsilon_{p}, \ud{l}_{p}}(f(\bar{u}_p,\cdot))$ and the corresponding $\bar{A}_p^{-1}$ is $\bar{\alpha}_p\bar{H}_p$
\STATE Update
$$
\bar{ u}_{p+1}=\bar{ u}_p+\bar{\alpha}_p\bar{ v}_p
$$
\STATE Compute $\bar{G}_{p+1}=[\mathcal{Q}_{\epsilon_{p+1}, \ud{l}_{p+1}}(D_{i}f(\bar{u}_{p+1},\cdot))]_{n\times 1}$
\STATE Define $\bar{s}_p:=\bar{ u}_{p+1}-\bar{ u}_p$ and $\bar{y}_p:=\bar{G}_{p+1}-\bar{G}_p$
\STATE Update
$$
\bar{H}_{p+1}=\left(I-\frac{\bar{s}_p\bar{y}_p^T}{\bar{s}_p^T\bar{y}_p}\right)\bar{H}_p\left(I-\frac{\bar{y}_p\bar{s}_p^T}{\bar{s}_p^T\bar{y}_p}\right)+\frac{\bar{s}_p\bar{s}_p^T}{\bar{s}_p^T\bar{y}_p}
$$
\STATE Set p:=p+1
\ENDWHILE
\STATE Output $\bar{ u}_p$ and $\bar{F}_p:=\mathcal{Q}_{\epsilon_{p}, \ud{l}_{p}}(f(\bar{u}_p,\cdot))$
\end{algorithmic}
\end{algorithm}

In order to illustrate our method, we provide the following 2D example. We will look into this example and it will also be used to explain the idea of next two sections.
\begin{example}\label{toy_exam}
We consider the following minimization problem
\begin{equation}
\min_{u\in \mathbb{R}} F(u)
\end{equation}
where $F(u)=\operatorname{\mathbb{E}}\left[u^2+(W_1^2+10W_2^2)u\right]$.
$W_1$ and $W_2$ are i.i.d. random variables. The integrand satisfies assumption $\ref{assumption}$. Moreover, the objective function is strictly convex in this example, so we conclude that there is a unique minimizer of this problem. By using the linearity of the expectation, the minimizer of the problem is
\begin{equation}
u^*=-\frac{\operatorname{\mathbb{E}} [W_1^2]+10\operatorname{\mathbb{E}} [W_2^2]}{2}.
\end{equation}

%The probability density function of Beta distribution is
%\begin{equation}
%p(w_i,\alpha,\beta)=\frac{w_i^{\alpha-1}(1-w_i)^{\beta-1}}{B(\alpha,\beta)}
%\end{equation}
%where $B$ is the beta function. The mean and variance of $Beta(\alpha,\beta)$ can be found in standard textbooks on probability. They are
%\begin{equation}
%\mathbb{E}W_i=\frac{\alpha}{\alpha+\beta},\ \ \ \mathrm{Var}\ W_i=\frac{\alpha\beta}{(\alpha+\beta)^2(\alpha+\beta+1)}.
%\end{equation}

\begin{figure}[h]
\centering
\subfigure[Solution error]{
\scalebox{0.35}[0.35]{\includegraphics{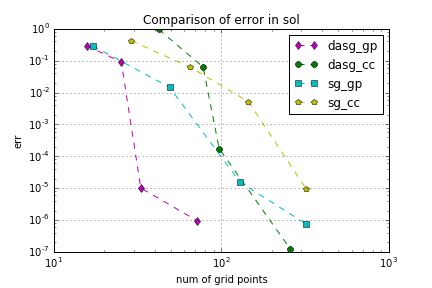}}
}
\subfigure[Solution error]{
\scalebox{0.35}[0.35]{\includegraphics{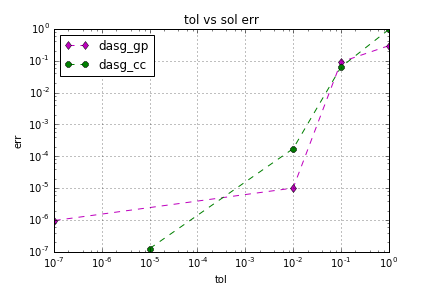}}
}
\subfigure[Objective error]{
\scalebox{0.35}[0.35]{\includegraphics{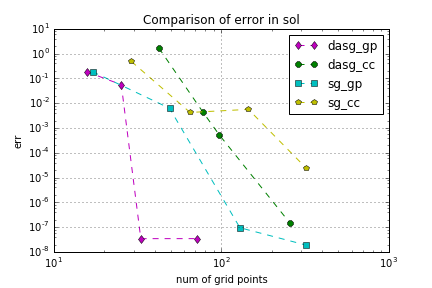}}
}
\subfigure[Objective error]{
\scalebox{0.35}[0.35]{\includegraphics{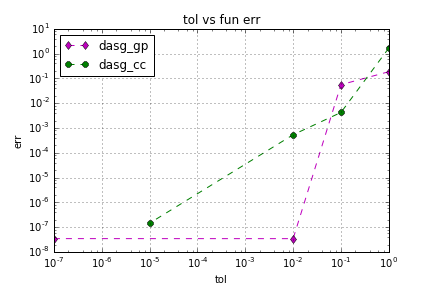}}
}

\caption{\footnotesize{Computational results for the 2D problem: (a) errors of computed minimizer vs. number of grid points used in each iteration(average for the dimension-adaptive sparse grid method).(b) errors of computed minimizer vs. $\epsilon$ in the termination condition of the dimension-adaptive algorithm. (c) errors of computed minimum vs. number of grid points used in each iteration. (d) errors of computed minimum vs. $\epsilon$ in the termination condition of the dimension-adaptive algorithm. We compare the dimension-adaptive sparse grid based on Gauss--Patterson(dasg$\_$gp) and Clenshaw--Curtis(dasg$\_$cc) with sparse grid based on Gauss--Patterson(sg$\_$gp) and Clenshaw--Curtis(sg$\_$cc).}}\label{toy_err}
\end{figure}
%We write down the explicit form of the gradient $G$ as follow
%\begin{equation}
%G(u)=\int_0^1\int_0^1 [2u+(w_1+w_2)]p(w_1,\alpha,\beta)p(w_2,\alpha,\beta)\ dw_1dw_2,
%\end{equation}
%where 
%\begin{equation}
%p(w_i,\alpha,\beta)=\frac{w_i^{\alpha-1}(1-w_i)^{\beta-1}}{B(\alpha,\beta)}
%\end{equation}
%is the probability density function of the Beta distribution. $B$ is the beta function here. 
In particular, here we further assume
\begin{equation}
W_k\sim\ \operatorname{Beta}(\alpha,\beta),\ k=1,2.
\end{equation}
with $\alpha=5,\beta=5$. The exact minimizer is then $u^*=-1.5$ and the minimum is $-2.25$.

In Figure $\ref{toy_err}$, we apply the Algorithm $\ref{bfgs_discretise}$ to solve the problem. We use forward difference to approximate the derivatives of the integrand. We compare the performance of the dimension-adaptive sparse grid quadrature and the sparse grid quadrature. For the OTDM based on the sparse grid quadrature, we fix the level $l$ for each run which results in the same choices of the quadrature rule $Q_p$ and $Q_{p,1}$ for any $p$. For the OTDM based on 
the dimension-adaptive sparse grid quadrature, we also fix the $\epsilon$ and $\ud{l}$ in the while condition $\ref{term}$. However, the choices of the quadrature rule are no longer the same for all $Q_p$ and $Q_{p,1}$. This is because the underlying downsets are not necessary to be the same for all $Q_p$ and $Q_{p,1}$. We see in Figure $\ref{toy_err}$ that the convergence rates are improved for both Clenshaw--Curtis and Gauss--Patterson when we apply the dimension-adaptive method. Especially, the average number of the grid points used in each iteration is substantially reduced for the same accuracy in the solution and the objective when Gauss-- Patterson is used as 1D rule. For Clenshaw--Curtis, we can also see this pattern, moreover, higher accuracy are obtained for both computed minimizer and minimum. 

\end{example}

\begin{figure}[h]\label{OTDM_DTOM}
\centering
\subfigure[Solution error]{
\scalebox{0.35}[0.35]{\includegraphics{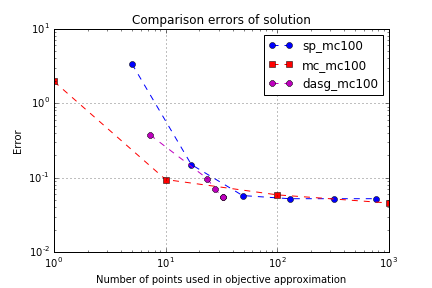}}
}
\subfigure[Objective error]{
\scalebox{0.35}[0.35]{\includegraphics{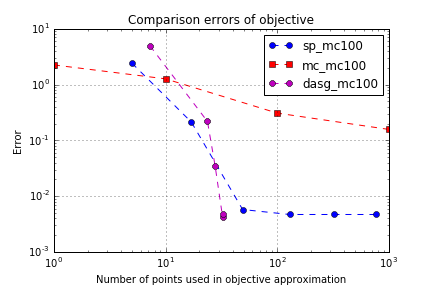}}
}
\subfigure[Solution error]{
\scalebox{0.35}[0.35]{\includegraphics{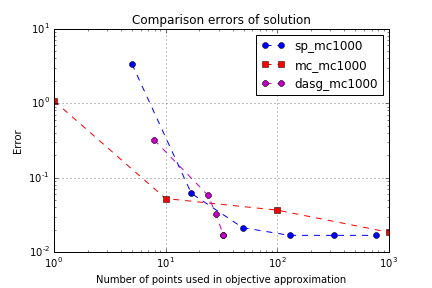}}
}
\subfigure[Objective error]{
\scalebox{0.35}[0.35]{\includegraphics{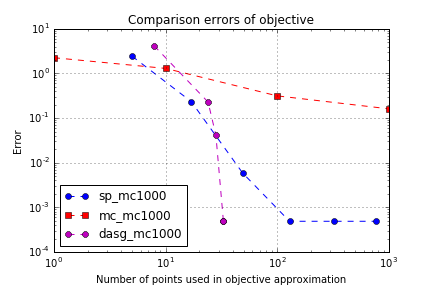}}
}

\caption{\footnotesize{Computational results for the toy problem. Here we use the exact expression of the gradient. Monte Carlo method is used to approximate the gradient and three different quadrature methods are used to approximate the objective function. (a) number of grid points vs. errors of computed minimizers(average of $10$ runs of each method with gradient approximated by $100$ points Monte Carlo method. (b)number of grid points vs. errors of computed minimum. (c)number of grid points vs. errors of computed minimizers($1000$ points Monte Carlo). (d)number of grid points vs. errors of computed minimum.}}
\end{figure}

The exact expression of the gradient in the previous 2D Example $\ref{toy_exam}$ is
\begin{equation}
\begin{aligned}
G(u)&=\nabla \operatorname{\mathbb{E}}[u^2+(W_1^2+10W_2^2)u]\\
&=\int_0^1\int_0^1 [2u+(w_1+10w_2)]p(w_1,\alpha,\beta)p(w_2,\alpha,\beta)\ dw_1dw_2
\end{aligned}
\end{equation}
where $p$ is the probability density function.
In Figure $\ref{OTDM_DTOM}$, we again solved problem with the BFGS method. We use $100$ points and $1000$ points Monte Carlo method to approximate the integral $G(u)$ respectively. The sparse grid, the dimensional-adaptive sparse grid based on Gauss--Patterson 1D quadrature and the Monte Carlo are used in approximating the objective function. We intentionally choose the same random points for Monte Carlo in objective approximation with those in gradient approximation when number of points equals to $100$ in (a), (b) and $1000$ in (c), (d). Thus, according to the propositions, the Monte Carlo surrogate methods are actually used. We can see from the Figure $\ref{OTDM_DTOM}$, both convergence performances of sparse grid and dimensional adaptive sparse grid are better than the Monte Carlo method (include those points which is actually surrogate method). The amount of work will be substantially reduced in the objective function if we apply the dimension-adaptive sparse grid in computing the objective. In addition, the dimension-adaptive sparse grid method performs better than sparse grid method in computing both solution and objective.

\section{Convergence analysis}
%If the problem is high dimensional in $\Omega$, the computational cost will be huge when we require a high accuracy approximation to the gradient $\nablaF(\ud{u_i})$. From the pseudo code, we need to compute the gradient in each iteration. If high accuracy is required in every iteration, the computational cost will be unaffordable. On the other hand, the convergence of Newton-type method may breakdown if we have bad approximation to the gradient, so we can't compute the gradient in very low accuracy. Thus, it is important for us to know the 'suitable' low accuracy of the approximation of the gradient such that the Newton-type method still works. Once we find out this 'suitable' accuracy we want, we can use $(\ref{bound})$ to decide the corresponding termination condition in our dimension adaptive sparse grid algorithm. 
We now study the sequence generated by the algorithm introduced in the previous section. Because of the presence of the errors (Both truncated errors and rounding errors), the actual sequences $\left\{\bar{u}_p\right\}$ generated by the Newton-type method is
\begin{equation}\label{perturb1}
\bar{u}_{p+1}=\left(I+E_0(\bar{u}_p)\right)\left(\bar{u}_p+\bar{\theta}(\bar{u}_p)\right),\ \ \ p=0,1,\dots,.
\end{equation}
In the above formula, $E_0(\bar{u}_p)$ is the relative error of the addition of $\bar{u}_p$ and $\bar{\theta}(\bar{u}_p)$. The norm of this error is of the order of the machine epsilon and thus very small, so this error is neglected in our analysis. $\bar{\theta}(\bar{u}_p)$ is the exact solution of the following linear system
\begin{equation}\label{az=g}
\bar{A}_pz=-\bar{G}_p
\end{equation}
where
\begin{equation}\label{perturb2}
\begin{aligned}
\bar{A}_p&=\hat{A}_p+E_1(\bar{u}_p)\\
\bar{G}_p&=G(\bar{u}_p)+E_2(\bar{u}_p).
\end{aligned}
\end{equation}
$\hat{A}_p$ is an approximation of the Hessian. Since $\hat{A}_p$ is computed by using the information of previously computed $\bar{u}_t$ and $\hat{A}_t$, $t<p$, we use the different notation $\hat{A}_p$ rather than $A_p$ here. $E_1(\bar{u}_p)$ and $E_2(\bar{u}_p)$ are errors occur when we approximate $\hat{A}_p$ and $G(\bar{u}_p)$. 

In \cite{YP}, the author shows the crucial condition for the convergence of such perturbed Newton-type method is 
\begin{equation}\label{condition}
\eta_p=\vert\vert\bar{A}_p^{-1}G'(\bar{u}_p)-I\vert\vert+\frac{\vert\vert\bar{A}_p^{-1}(\bar{G}_p-G_p)\vert\vert}{\vert\vert(G^{'}(\bar{u}_p))^{-1}G_p\vert\vert}<1
\end{equation}
The following theorem~\cite{YP} provides more details on the convergence result
\begin{thm}\label{convergence}
If $\vert\vert E_0\vert\vert=0$ and $\eta_p\le \eta<1$ for all $p=0,1,\dots,$ and if $\bar{u}_0$ satisfies
\begin{equation}
\vert\vert\bar{u}_0-u^*\vert\vert<\frac{2(1-\eta)}{(3-\eta)\mu},
\end{equation}
then the perturbed Newton-type method $(\ref{perturb1})-(\ref{perturb2})$ produces a sequence $\left\{\bar{u}_p\right\}$ which converges to $u^*$.
\end{thm}
Because of the presence of the term $\vert\vert(G^{'}(\bar{u}_p))^{-1}G_p\vert\vert$ in the denominator, we can expect that the inequality $(\ref{condition})$ will finally be violated and the convergence theory will fail. Therefore, there exists some neighborhood of $u^*$ such that the sequence generated by the Newton-type method converges outside this neighborhood while the behavior of the sequence generated by the successive iterations is not predictable if they are inside the neighborhood. This also suggests further iterations will not improve the accuracy and we can stop the algorithm once the  $\bar{u}_p$ reaches the neighborhood. 

In order to find out when to stop further iterations, we look into the two terms in $(\ref{condition})$. When $G'(\bar{u}_p)$ can be well approximated and $(\ref{az=g})$ can be solved with relative high accuracy, the first term can be kept small. This means the errors in the Jacobian approximation are tolerable if
\begin{equation}\label{part_1_condition}
\tau_p=\vert\vert\bar{A}_p^{-1}G'(\bar{u}_p)-I\vert\vert\ll 1.
\end{equation} 
%In this case, the second term can be further written as
%\begin{equation}\label{upp bound}
%\frac{\vert\vert G'(\bar{u}_i)^{-1}(\bar{G}_i-G_i)\vert\vert}{\vert\vert(G^{'}(\bar{u}_i))^{-1}G_i\vert\vert}\leq\kappa(G'(\bar{u}_i))\frac{\vert\vert G_i-\bar{G}_i\vert\vert}{\vert\vert G_i\vert\vert},
%\end{equation}
%where $\kappa(G'(\bar{u_i}))$ is the condition number of $G'(\bar{u_i})$. Therefore, the Newton-type method may not be convergent when
%\begin{equation}
%\kappa(G'(\bar{u}_i))\frac{\vert\vert G_i-\bar{G}_i\vert\vert}{\vert\vert G_i\vert\vert}>1.
%\end{equation} 
%When $U$ is a subspace of $\mathbb{R}$, the condition number $\kappa(G')=1$. In this case, we can roughly estimate the smallest number of iterations $i_{s}$ for $\left\{\bar{u_i}\right\}$ to reach the neighborhood we mentioned before by solving $\vert\vert G_i-\bar{G}_i\vert\vert>\vert\vert G_i\vert\vert$. When $U$ is a high dimensional space, the condition number may not equal to $1$. We need to estimate the condition number first. By using the perturbation lemma[ref], we have 
%\begin{equation}
%\kappa(G'(\bar{u}_i))\le \kappa(\bar{A_i})\frac{1+\tau_i}{1-\tau_i}\approx \kappa(\bar{A_i}).
%\end{equation}
%One can get an estimation of $\kappa(\bar{A}_i)$ when solving the system equation $(\ref{az=g})$, we denote it as $\tilde{\kappa}(\bar{A}_i)$. Then the estimate of $i_s$ can be followed by solving inequality $\tilde{\kappa}(\bar{A}_i)\vert\vert G_i-\bar{G}_i\vert\vert>\vert\vert G_i\vert\vert$. We will stop the iteration when $i>i_s$ since further iterations will not improve the computational result. 
\begin{cor}~\cite{YP}\label{breakdown}
If $(\ref{part_1_condition})$ holds, then the convergence of the perturbed Newton-type method might breakdown when
\begin{equation}\label{YPcondition}
\lVert\bar{G}_p-G_p\rVert>\frac{\lVert \bar{G}_p\rVert}{1+\kappa(\bar{A}_p)}.
\end{equation}
where $\kappa(\bar{A}_p)$ is the condition number of $\bar{A}_p$.
\end{cor}

%The magnitude of the $i_s$ is dependent on the problem itself and also how accurate the integrals and derivatives are approximated. For the same problem, if these integrals and derivatives are approximated in relative high accuracy, we will get larger $i_s$ and the Newton-type method will finally achieve a more accurate minimizer. On the other hand, if we can only achieve an low accuracy approximation for these integrals and derivatives which might happen during solving very high dimensional problems, the associated number of iterations in Newton-type method should be small and the final solution can only be solved in low accuracy. We further explore the relation between $i_s$ and the accuracy of approximation to the integrals and derivatives by doing some analysis.

According to the Corollary $\ref{breakdown}$, the accuracy of computed minimizer depends on how accurate the gradient can be approximated. For the same problem, if the gradients can be approximated in higher accuracy, which means we can get smaller $\vert\vert \bar{G}_p-G_p\vert\vert$, then it is likely that more effective iterations can be carried on in the Newton-type method according to $(\ref{YPcondition})$. Thus, we can obtain a better solution. The problem is it will be very expansive to get a moderately accurate approximation to the gradient if we are dealing with high dimensional problems. Thus, there is no need to do iteration $(\ref{iter})$ many times for this kind of problems and $(\ref{YPcondition})$ actually provide us with a simple stopping criterion. In the stopping criterion $(\ref{YPcondition})$, the condition number $\kappa(\bar{A}_p)$ can be computed during the iteration with little cost. However, it will be difficult to obtain the exact value $\vert\vert \bar{G}_p-G_p\vert\vert$ and $\vert\vert G_p\vert\vert$. Following theorem gives an estimation of the ratio of
$\vert\vert \bar{G}_p-G_p\vert\vert$ to $\vert\vert G_p\vert\vert$ and  provide more detailed stopping criterion for our method.

\begin{thm}\label{stopping criterion}
Consider solving the stochastic optimisation problem $(\ref{F})$, we assume here the integrand satisfies
\begin{itemize}
\item[(i)] $f(\cdot,W)\in C^2(U)$ for all $W\in \Omega$.
\item[(ii)] $f(u,\cdot)\in H^r(\Omega)(r\geq 1)$ for all $u\in U$.
\end{itemize}
We apply the dimension-adaptive sparse grid($\epsilon$) to approximate the integrals and forward difference method to approximate the derivatives in the computation. We further assume all the components of the gradient at $p$th iteration are computed by the same downset as that in computing the objective function at $p$th iteration. Then, we have
\begin{equation}\label{error_bound_g}
\lVert\bar{G}_p-G_p\rVert\le \lVert\mathcal{E}_p^1\rVert+\lVert\mathcal{E}_p^2\rVert
\end{equation} 
where the $q$th element of $\mathcal{E}_p^1$ and $\mathcal{E}_p^1$ are
\begin{equation}
\begin{aligned}
\mathcal{E}_{p,q}^1&=K_{p,q}^1h\\
\mathcal{E}_{p,q}^1&=K_{p,q}^2 2^{-l_pr}+K_{p,q}^3\sum_{\ud{i}\in \ud{L}_{p}\backslash\ud{I}_p}2^{-\vert\ud{i}\vert r}
\end{aligned}
\end{equation}
where $K_{p,q}^i,\ i=1,2,3$ are constants. Therefore, our method based on dimension-adaptive sparse grid method might breakdown when
\begin{equation}\label{breakdown_con}
\lVert\mathcal{E}_p^1\rVert+\lVert\mathcal{E}_p^2\rVert>\frac{\lVert \bar{G}_p\rVert}{1+\kappa(\bar{A}_p)}.
\end{equation}
\end{thm}
\begin{proof}
We use forward difference quotient
\begin{equation}
\tilde{G}_{p,q}=\frac{F(\bar{u}_p+h\ud{e}_q)-F(\bar{u}_p)}{h}
\end{equation}
to approximate the $q$th component $G_{p,q}$ of $G_p$, where $h$ is the difference increment. We further denote the $q$th component of $\bar{G}_p$ as $\bar{G}_{p,q}$ which can be expressed as
\begin{equation}
\bar{G}_{p,q}=\frac{\bar{F}(\bar{u}_p+h\ud{e}_q)-\bar{F}(\bar{u}_p)}{h}
\end{equation}
where $\bar{F}(\bar{u}_p)=Q_{\ud{I}_p}(f(\bar{u}_p,\cdot))$.

By using triangle inequality, we have
\begin{equation}\label{triangle_in}
\lVert \bar{G}_p-G_p\rVert\le \lVert \bar{G}_p-\tilde{G}_p+\tilde{G}_p-G_p\rVert\le \lVert\bar{G}_p-\tilde{G}_p\rVert+\lVert\tilde{G}_p-G_p\rVert.
\end{equation}
Since $f(\cdot,W)\in C^2(U)$ for any $W\in \Omega$, we can get $F\in C^2(U)$. Thus, for the $q$th component of $\tilde{G}_p-G_p$, we have
\begin{equation}\label{second_term}
\lvert \tilde{G}_{p,q}-G_{p,q}\rvert\le K_{p,q}^1h.
\end{equation}
For the term $\lVert\bar{G}_p-\tilde{G}_p\rVert$, we have
\begin{equation}
\begin{aligned}
&\bar{G}_{p,q}-\tilde{G}_{p,q}\\
=&\frac{1}{h}(\bar{F}(\bar{u}_p+h\ud{e}_q)-\bar{F}(\bar{u}_p))-\frac{1}{h}(F(\bar{u}_p+h\ud{e}_q)-F(\bar{u}_p))\\
=&\frac{1}{h}(Q_{\ud{I}_p}(f(\bar{u}_p+h\bar{e}_q,\cdot))-Q_{\ud{I}_p}(f(\bar{u}_p,\cdot)))-\frac{1}{h}(I(f(\bar{u}_p+h\bar{e}_q,\cdot))-I(f(\bar{u}_p,\cdot)))\\
=&Q_{\ud{I}_p}\left(\frac{1}{h}(f(\bar{u}_p+h\bar{e}_q,\cdot)-f(\bar{u}_p,\cdot))\right)-I\left(\frac{1}{h}(f(\bar{u}_p+h\bar{e}_q,\cdot)-f(\bar{u}_p,\cdot))\right)\\
=&(Q_{\ud{I}_p}-I)\left(\frac{1}{h}(f(\bar{u}_p+h\bar{e}_q,\cdot)-f(\bar{u}_p,\cdot))\right).
\end{aligned}
\end{equation}
We used the assumption of the same downset $I_p$ in the second equality. The third equality is due to the linearity of the operator $Q_{\ud{I}_p}$ and $I$. Since $f(u,\cdot)\in H^r(\Omega)$ for any $u\in U$, we know that the function
\begin{equation}
\frac{1}{h}(f(\bar{u}_p+h\bar{e}_q,\cdot)-f(\bar{u}_p,\cdot))\in H^r(\Omega).
\end{equation}
Thus, we get following upper bound
\begin{equation}\label{fisrt_term}
\lvert \bar{G}_{p,q}-\tilde{G}_{p,q}\rvert\le K_{p,q}^2 2^{-l_pr}+K_{p,q}^3\sum_{\ud{i}\in \ud{L}_{p}\backslash\ud{I}_p}2^{-\vert\ud{i}\vert r}
\end{equation}
Combining inequalities $(\ref{triangle_in})$, $(\ref{fisrt_term})$ and $(\ref{second_term})$, we get $(\ref{error_bound_g})$.

By using Corollary $\ref{breakdown}$, we get the breakdown condition $(\ref{breakdown_con})$.
\end{proof}

\begin{rmk}
Theorem $(\ref{stopping criterion})$ actually provides us with rough stopping criterion. It can be checked at each iteration in Newton-type method if one can get reasonable estimations of the constants $K_{p,q}^1$, $K_{p,q}^2$, $K_{p,q}^3$. The $h$, $l_p$, $\ud{L_p}$ and $\ud{I_p}$ can be obtained during the computation. The second terms of $\mathcal{E}_{p,q}^2$ can also be replaced by the bound in corollary $(\ref{fin_err_bound})$. 
The advantage of doing this is that instead of estimating the constant $K_{p,q}^3$, we can get a reasonable stop by tuning $\rho$.
\end{rmk}

\section{Stopping Criterion}
In this section, we discuss the stopping criterion for our method. It is important to know when we should stop the Newton-type method for fixed tolerance $\epsilon$ in the termination condition of the dimension-adaptive sparse grid method. A good stopping criterion can save a lot of computational cost. This is because if the problem is high dimensional, each iteration will be very expansive to compute even if we apply the dimension-adaptive sparse grid method to reduce the cost in computing the related integrals. Another reason is it is possible that further iterations can not improve the accuracy of the solution. It might happen that the computed solution becomes even worse after we increase more iterations. Also, computation of the termination condition of some existing optimization solvers can involve a great number of evaluations of objective functions and gradients, such as the Wolfe's rule used in BFGS method in scipy.optimize package. Sometimes, computing the termination condition is even expansive than computing the solution itself. The study of the stopping criterion can also shed light on how to choose the tolerance  
$\epsilon$ in the termination condition of the dimension-adaptive sparse grid method for a specific problem with given requirement on the accuracy of the solution.

Though Theorem $(\ref{stopping criterion})$ gives us some hints to derive a stopping criterion, accurate estimations of the coefficients are required. However, this is usually difficult especially for high dimensional problems. In addition, the estimation method can vary for different problems, thus it will be complex to obtain a general stopping criterion from the Theorem. 

Here we provide another way to find the stopping criterion. It works better and much easier to implement in practice. Suppose $\bar{u}_p$ is the approximated minimizer generated by some Newton-type methods after $p$th iteration. When $\bar{u}_p$ is close enough to the exact minimizer $u^*$, we have following Taylor expansion
\begin{equation}
F(\bar{u}_p)=F(u^*)+\nabla F(u^*)(\bar{u}_p-u^*)+(\bar{u}_p-u^*)^T\nabla^2 F(u^*)(\bar{u}_p-u^*)+o(\Vert\bar{u}_p-u^*\Vert^2).
\end{equation}
Since $u^*$ is the minimizer of $F$, we have $\nabla F(u^*)=0$. $F$ is convex, therefore $\nabla^2 F(u^*)$ is positive semidefinite. Thus, we have $\Vert\bar{u}_p-u^*\Vert^2$ increases(decreases) with $i$ as $F(\bar{u}_p)$ increases(decreases). However, the exact value of $F(\bar{u}_p)$ is usually not easy to get. Therefore, in our approach, instead of using the exact function value at $\bar{u}_p$, we use the value  of some high accuracy approximations of the function. If we denote the high accuracy approximation by dimension-adaptive sparse grid quadrature with $\epsilon$ in its termination condition at $u$ as $F_{\epsilon}(u)=Q_{\ud{I}}f(u,\cdot)$, then we can decide when to stop the Newton-type method by studying the trend of $F_{\epsilon}(\bar{u}_p)$.

The advantage of this stopping criterion is the computation does not  involve estimation of the error $\Vert G-\bar{G}\Vert$. We only need to compute $F_{\epsilon}(\bar{u}_p)$ for some smaller $\epsilon$s with the same algorithm which used to compute the approximation of the objective function. Also, the additional computational cost for computing $F_{\epsilon}(\bar{u}_p)$ is affordable in most cases for even high dimensional problems. This is because we use Newton-type method as our optimisation solver, so the number of iterations will not be too large. Also, the computational cost of getting such stopping criterion is much lower than that of computing the gradient in each iteration when the dimension of $U$ is high.

\begin{figure} 
\centering 
\subfigure[Gauss Patterson]{ 
\begin{minipage}[b]{0.30\linewidth} \includegraphics[width=1\linewidth]{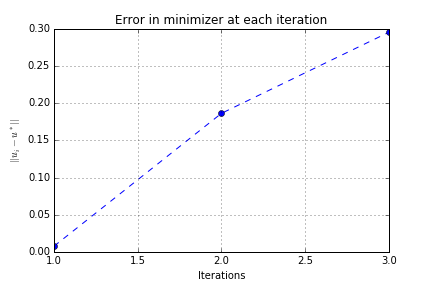}\vspace{4pt} \includegraphics[width=1\linewidth]{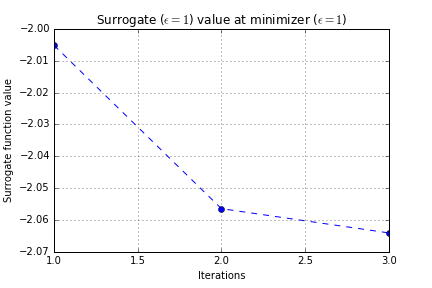}\vspace{4pt} \includegraphics[width=1\linewidth]{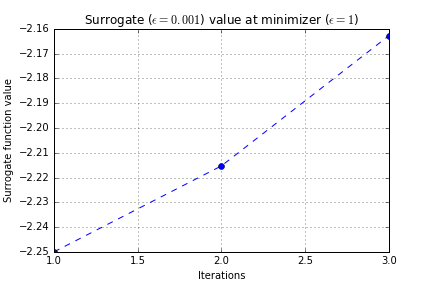}\vspace{4pt} \includegraphics[width=1\linewidth]{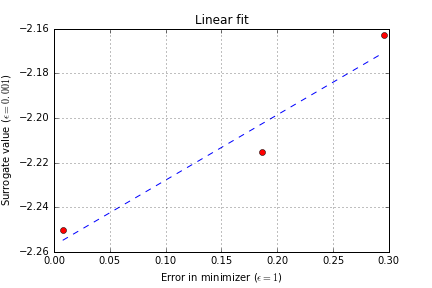} 
\end{minipage}} 
\subfigure[Clenshaw Curtis]{ 
\begin{minipage}[b]{0.30\linewidth} \includegraphics[width=1\linewidth]{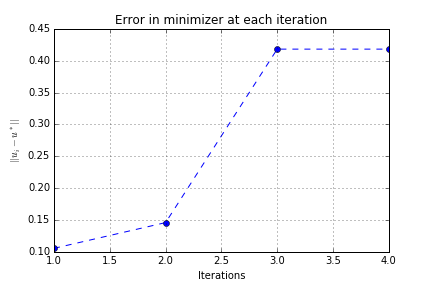}\vspace{4pt} \includegraphics[width=1\linewidth]{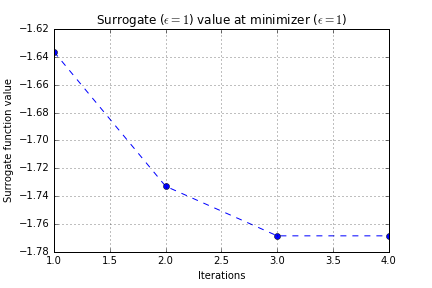}\vspace{4pt} \includegraphics[width=1\linewidth]{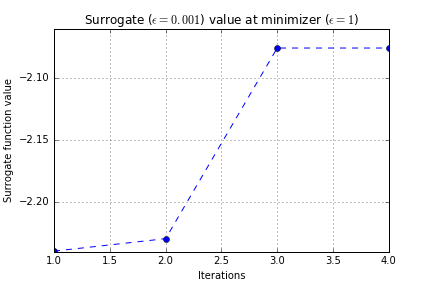}\vspace{4pt} \includegraphics[width=1\linewidth]{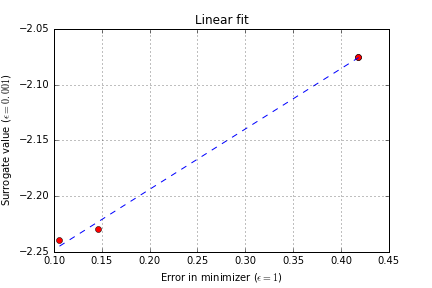} 
\end{minipage}} 
\subfigure[Trapezoidal Rule]{ 
\begin{minipage}[b]{0.30\linewidth} \includegraphics[width=1\linewidth]{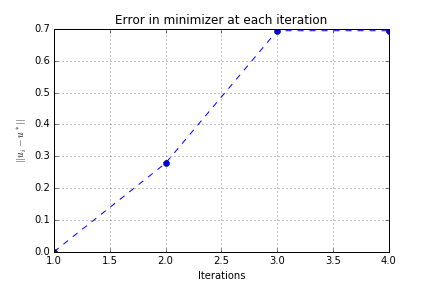}\vspace{4pt} \includegraphics[width=1\linewidth]{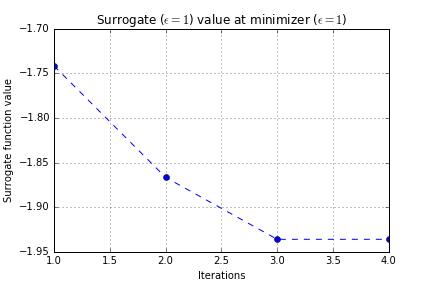}\vspace{4pt} \includegraphics[width=1\linewidth]{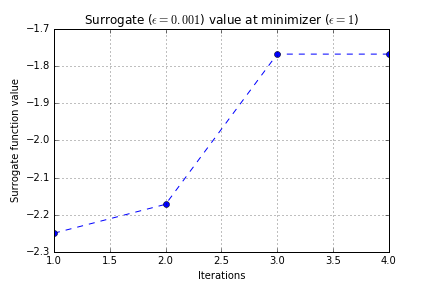}\vspace{4pt} \includegraphics[width=1\linewidth]{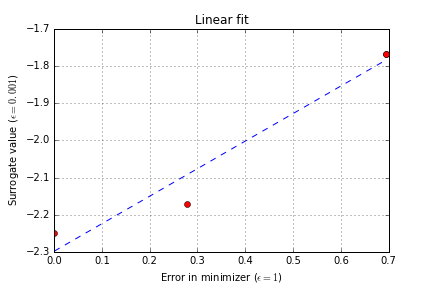} 
\end{minipage}} 
\caption{Solve the toy problem with $\epsilon=1$.}\label{stopping_1} 
\end{figure}

\begin{figure}[h] 
\centering 
\subfigure[Gauss Patterson]{ 
\begin{minipage}[b]{0.30\linewidth} \includegraphics[width=1\linewidth]{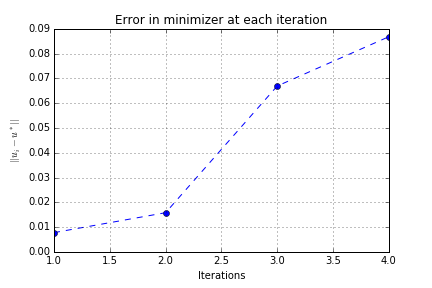}\vspace{4pt} \includegraphics[width=1\linewidth]{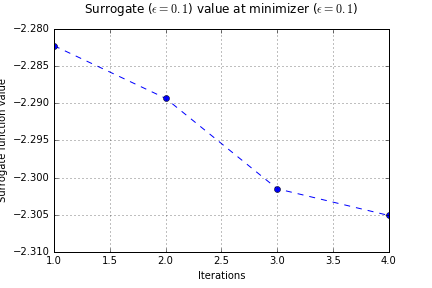}\vspace{4pt} \includegraphics[width=1\linewidth]{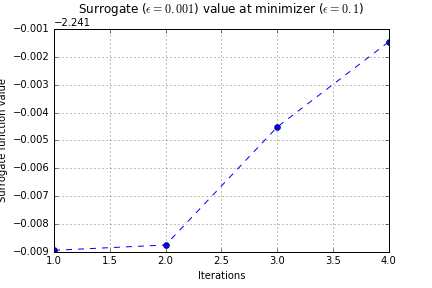}\vspace{4pt} \includegraphics[width=1\linewidth]{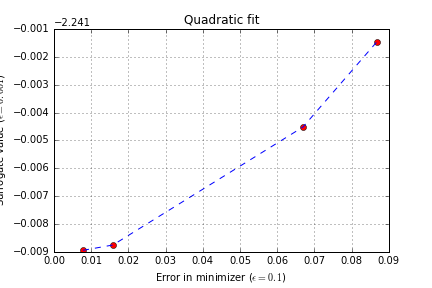} 
\end{minipage}} 
\subfigure[Clenshaw Curtis]{ 
\begin{minipage}[b]{0.30\linewidth} \includegraphics[width=1\linewidth]{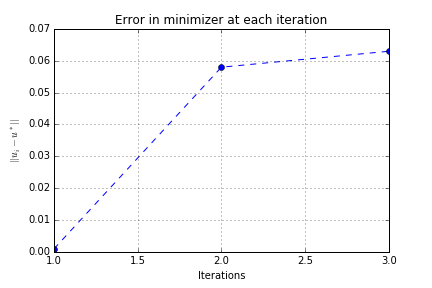}\vspace{4pt} \includegraphics[width=1\linewidth]{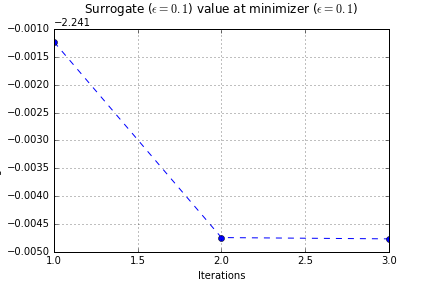}\vspace{4pt} \includegraphics[width=1\linewidth]{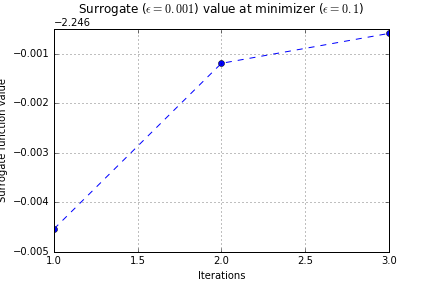}\vspace{4pt} \includegraphics[width=1\linewidth]{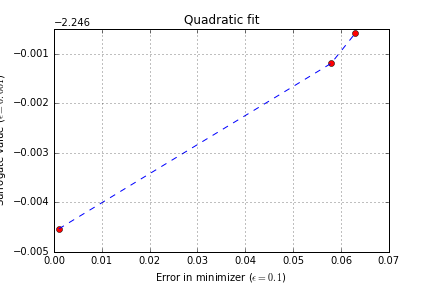} 
\end{minipage}} 
\subfigure[Trapezoidal Rule]{ 
\begin{minipage}[b]{0.30\linewidth} \includegraphics[width=1\linewidth]{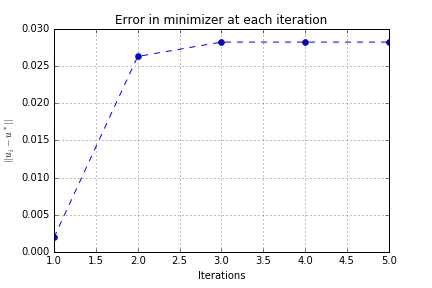}\vspace{4pt} \includegraphics[width=1\linewidth]{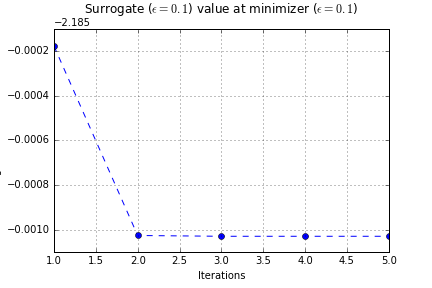}\vspace{4pt} \includegraphics[width=1\linewidth]{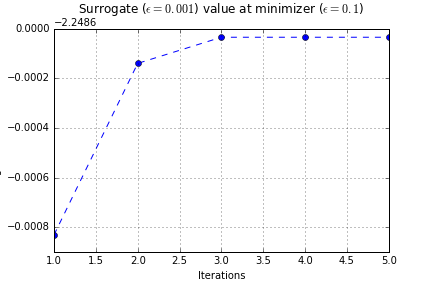}\vspace{4pt} \includegraphics[width=1\linewidth]{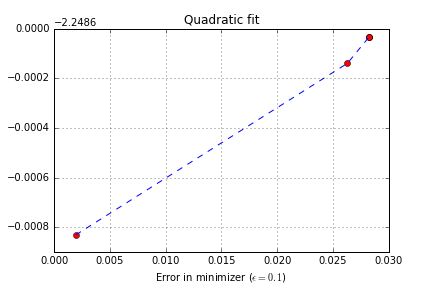} 
\end{minipage}} 
\caption{Solve the toy problem with $\epsilon=0.1$.}\label{stopping_2} 
\end{figure}

In fig \ref{stopping_1} and fig \ref{stopping_2}, we solve the toy problem by using surrogate method with $\epsilon=1$ and $\epsilon=0.1$ respectively. For both two figures, the subfigures in the first row show the relation between error $\vert u^*-\bar{u}_p\vert$ versus the number of iterations. We can see from these figures, we should stop the algorithm after first iteration for all three cases. If we further increase the iteration, the errors will grows. The reason for this is the gradients are approximated with relatively low accuracy. From the convergence theory, we know the convergence of the Newton-type method might break down in this case. The subfigures in the second row presents function values on the surrogate versus the number of iterations. As we expected, $F_{\epsilon}(\bar{u}_p)$ is decreasing. The subfigures in the third row show the function values of the surrogate functions with $\epsilon=0.001$ on $\bar{u}_p$. Compare the subfigures in the first row with the corresponding subfigures in the third row, we can see that the trends of the functions are the same for all three quadrature rules. Thus, we can predict when to stop the Newton-type method by studying the trends of the functions in the third row respectively. Our method successfully predicts the stopping times(after first iteration) for all three quadrature rules in this example. The subfigures in the fourth row, we fit the data points $(\vert u^*-\bar{u}_p\vert, F_{\epsilon/4}(\bar{u}_p))$ to a quadratic function. We find that the quadratic curves almost go through all the data points which suggests our error model is reasonable.

%\begin{figure} 
%\centering 
%\subfigure[Gauss Patterson]{ 
%\begin{minipage}[b]{0.30\linewidth} \includegraphics[width=1\linewidth]{gp_001_0.png}\vspace{4pt} \includegraphics[width=1\linewidth]{gp_001_1.png}\vspace{4pt} \includegraphics[width=1\linewidth]{gp_001_2.png}\vspace{4pt} \includegraphics[width=1\linewidth]{gp_001_3.png} 
%\end{minipage}} 
%\subfigure[Clenshaw Curtis]{ 
%\begin{minipage}[b]{0.30\linewidth} \includegraphics[width=1\linewidth]{cc_001_0.png}\vspace{4pt} \includegraphics[width=1\linewidth]{cc_001_1.png}\vspace{4pt} \includegraphics[width=1\linewidth]{cc_001_2.png}\vspace{4pt} \includegraphics[width=1\linewidth]{cc_001_3.png} 
%\end{minipage}} 
%\subfigure[Trapezoidal Rule]{ 
%\begin{minipage}[b]{0.30\linewidth} \includegraphics[width=1\linewidth]{tra_001_0.png}\vspace{4pt} \includegraphics[width=1\linewidth]{tra_001_1.png}\vspace{4pt} \includegraphics[width=1\linewidth]{tra_001_2.png}\vspace{4pt} \includegraphics[width=1\linewidth]{tra_001_3.png} 
%\end{minipage}} 
%\caption{Solve the toy problem with $\epsilon=0.01$.} 
%\end{figure}

\section{High dimensional examples}
\subsection{50D additive example}
Consider following minimization problem
\begin{example}
\begin{equation}
\min_{u\in U} \displaystyle \mathop{\mathbb{E}}\left[\sum_{i=1}^d \exp{(-u_iW_i^2)}\right]
\end{equation}
where $W_i$ are i.i.d random variables which subject to uniform distribution $U(0,1)$ and the set $U=[0,1]^d$. Thus the integral form of the objective function can be written as
\begin{equation}
F(u)=\int_{[0,1]^d}\sum_{i=1}^d e^{-u_iw_i^2}\ dw.
\end{equation}
The gradient $G(u)$ is
\begin{equation}
\begin{aligned}
G(u)&=\nabla F(u)\\
&=-\left[\int_{[0,1]^d}w_1^2e^{-u_1w_1^2}\ dw_1,\dots,\int_{[0,1]^d}w_d^2e^{-u_dw_d^2}\ dw_d\right],
\end{aligned}
\end{equation}
so we have $G(u)\le \ud{0}$ for any $u\in [0,1]^d$ and thus the exact minimizer is $u^*=(1,\dots,1)$ for this problem.

The reference objective function value can be computed by 
\begin{equation}
F(u)=\int_{[0,1]^d}\sum_{i=1}^d e^{-u_iw_i^2}\ dW=\sum_{i=1}^d\int_{[0,1]} e^{-u_iw_i^2}\ dw_i.
\end{equation}
At the minimizer, we have the exact objective
\begin{equation}
F(u^*)=d\int_{[0,1]}e^{-w_i^2}\ dw_i=d(\mathcal{F}(1)-\mathcal{F}(0))
\end{equation}
where $\mathcal{F}$ is the cumulative distribution function.

\begin{figure}[h]
\centering
\begin{tikzpicture}
\begin{loglogaxis}[
width=0.7\textwidth,
height=0.5\textwidth,
title={Error in objective function},
xlabel={Average number of quadrature points},
ylabel={Error},
ymajorgrids=true,
grid style=dashed,
legend pos=south west
]
\addplot[dashed,blue,mark=*,mark options={solid,fill=blue}]
    coordinates {
    (1,0.979726590157) (10,0.2520708608) (100,0.1574037672) (1000,0.0323470791237) (10000,0.0122928437404) (100000,0.00348271281182)
};
\addplot[dashed,red,mark=square*,mark options={solid,fill=red}]
coordinates {    (1.0,1.59883215027) (101.0,0.0178047229758) (5364.0,2.29230634474e-05) (13492.0,5.63280864441e-06)
};
\addplot[dashed,brown,mark=pentagon*,mark options={solid,fill=brown}]
    coordinates {
    (1.0,1.59883215027) (101.0,0.000468074603418) (5101.0,7.73249055896e-07)
};
\legend{mc,dasgcc,dasggp}

\end{loglogaxis}
\end{tikzpicture}
\caption{\footnotesize {Compute the additive example with $d=50$.}}
\end{figure}
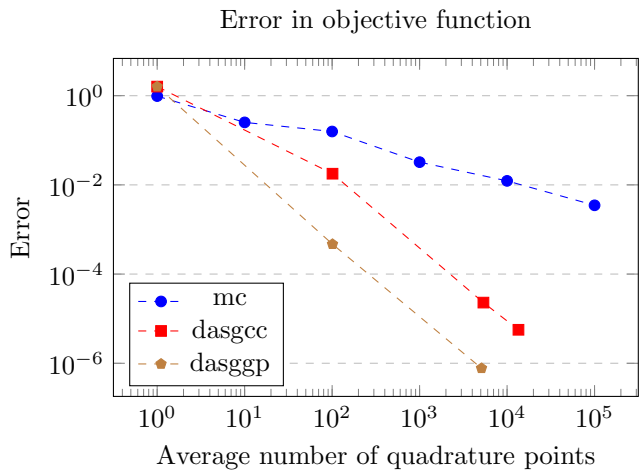

Here, we apply the 'optimize then discretise' approach to solve the problem. We use 'L-BFGS-B' method in scipy.optimize package as our solver. We apply the dimension-adaptive sparse grid method to computing the objective function while use the Monte-Carlo method to approximate the high dimensional integrals involved in computing the gradient.

We still have the approximated gradient $\bar{G}(u)\le 0$ for any $u\in [0,1]^d$. This is because the integrand of each entry of the gradient is non-positive and only positive weights are used in the Monte Carlo method.  Thus, no matter how many samples used in the Monte Carlo method, we will always get a descent direction at each step during optimization process. It is noteworthy here that both low level sparse grid method and dimensional-adaptive sparse grid method with large $\epsilon$ may change the sign of integral approximated and therefore lead to wrong search direction. 

In the example, we only use the Monte Carlo method with 10 samples to compute the gradient components. In order to increase the accuracy in minimum, when we approximate the objective function, we increase samples in the Monte Carlo method and decrease $\epsilon$ in the termination condition of the dimension-adaptive sparse grid. The result shows all of three methods achieve the exact minimizer as we expected. The errors in minimum of two dimension-adaptive approaches drop much faster than the Monte Carlo method.
\end{example}

%\subsection{Stochastic control}
\subsection{Application to stochastic control}
In this section, we illustrate our dimension-adaptive sparse grid method with an instance of a discrete time open-loop stochastic control problem. The general form of such control problem can be found in Bertsekas~\cite{Di}. The control problem is described by following discrete time dynamic system 
\begin{equation}\label{dynamic system}
x_{i+1}=\psi_i(x_i,u_i,w_i),\ i=0,\dots,d-1\\
\end{equation}
Here $x_i$ and $u_i$ are states and controls respectively where $x_0$ is given. $w_i$ are disturbances. Here we only consider a special case when the states, the controls and the disturbances are in one dimensional space. When the disturbances in the system are unknown, we  usually model them as i.i.d. random variables $W_i$ with given probability density function. In this case, the open-loop means the controls $u_i$ do not depend on the disturbances~\cite{BO} and we can further write the dynamic system in its random form:
\begin{equation}\label{dynamic system}
X_{i+1}=\psi_i(X_i,u_i,W_i),\ i=0,\dots,d-1.
\end{equation}

If we further define the vectors of states,controls and noises, i.e.,
\begin{equation}
X=(x_0,\dots,X_{d-1}),\ u=(u_0,\dots,u_{d-1}),\ W=(W_0,\dots,W_{d-1}),
\end{equation}
then we can rewrite the dynamic system as
\begin{equation}\label{Psi}
X=\Psi(X,u,W).
\end{equation}
where $\Psi$ is a function can be derived from $\psi_i$.
%together with $x_{d}=\psi_{d-1}(x_{d-1},u_{d-1},w_{d-1})$. 

Our task now is to determine what is the 'best' control for the dynamic system $(\ref{dynamic system})$ or $(\ref{Psi})$ to minimize the expected cost 
\begin{equation}
\mathbb{E}\left[\Phi(u,X)\right]
\end{equation}
where $\Phi$ is a given function. 

%In this section, we study a kind of discrete time open-loop stochastic control problem. The general form of a discrete time open-loop stochastic control problem is as follow. Minimize the objective function(the cost function) over $U$(the set contains all possible decisions), i.e.
%\begin{equation}
%\min_{u\in U} \mathbb{E}\left[\phi(u,X)\right],
%\end{equation}
%with constraints(dynamic system)
%\begin{equation}\label{ex_constraints}
%e(u,X(W),W)=0,\ \forall W\in \Omega.
%\end{equation}
%Here the $i$th component of $X$ and $u$ is the $i$th state variable and control variable at time step $i$ respectively. The randomness comes from a random vector $W$ defined on the probability space $(\Omega,\mathcal{B},\mathbb{P})$. In most case, the components of $W$ are i.i.d. continuous random variables with given probability density function $p(w_i)$. $\phi$ is a given function. Open-loop here means the control variables $u_i$ do not depend on $W$.

Here we focus on the case when $X$ can be solved explicitly from the dynamic system $(\ref{Psi})$, that is,
\begin{equation}
X=\xi(u,W).
\end{equation}
In this case, the original problem can be reduced into the standard form of the stochastic optimization problem discussed in the paper, namely,
\begin{equation}
\min_{u\in U}\mathbb{E}\left[h(u,W)\right],
\end{equation}
where
\begin{equation}
h(u,W)=\Phi(u,\xi(u,W)).
\end{equation}
The integral form of the expected cost and its surrogate with $N$ quadrature points are
\begin{equation}
\int_{\mathbb{R}^d}\Phi(u,\xi(u,w))p(w)\ dW\approx\sum_{j=1}^Nc_i\Phi(u,\xi(u,w_j))p(w_j).
\end{equation}
In order to illustrate the computational performance of our approach, we consider a classical example with linear dynamic system 
\begin{equation}\label{linear constraints}
X=AX+Bu+CW+x_0\ud{e_0}.
\end{equation}
and the quadratic objective function $\Phi$
\begin{equation}\label{quad_obj_fun}
\Phi(u,X)=u^TPu+x^TQx
\end{equation}
where $A$, $B$, $C$, $P$ and $Q$ are given $d\times d$ matrices and $x_0$ is the given initial value. By solving $(\ref{linear constraints})$, we get $\xi(u,W)=(I-A)^{-1}(Bu+CW+x_0\ud{e_0})$. Combining the expression of $\xi(u,W)$ with $(\ref{quad_obj_fun})$, we know that $h(u,W)$ is again a quadratic function.
%\begin{equation}
%u^TPu+\left[(I-A)^{-1}(Bu+CW+x_0\ud{e_0})\right]^{T}Q\left[(I-A)^{-1}(Bu+CW+x_0\ud{e_0})\right].
%\end{equation}

The exact solution can be derived by using the certainty equivalence principle~\cite{Di}. According to the principle, the solution of the stochastic control problem is the same as that of a corresponding deterministic problem when the objective function is quadratic and the constraints are linear. That means we can get the reference solution by numerically solving the deterministic problem(see appendix).

%\begin{figure}[h]\label{open loop problem}
%\centering
%
%\scalebox{0.5}[0.5]{\includegraphics{open_loop_minimizer_error_comparison.png}}
%
%
%\caption{\footnotesize{Computational results for the open-loop stochastic control problem with parameters $d=7$, $p_i=q_i=\frac{1}{d},\ \ \ i=0,\dots,d-1$, $a_i=(1+\frac{1}{d}),\ b_i=1,\ c_i=1,\ \ \ i=0,\dots,d-1$,$W_i\sim beta(2,3),\ \ \ i=0,\dots,d-1$.}}
%\end{figure}

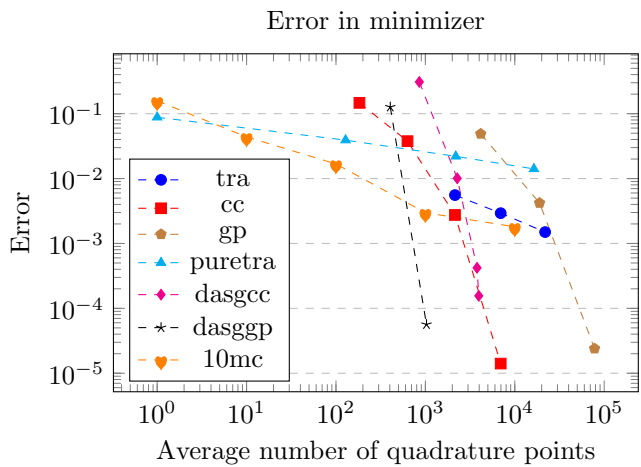
\begin{figure}[h]\label{open loop problem}
\centering
\begin{tikzpicture}
\begin{loglogaxis}[
width=0.7\textwidth,
height=0.5\textwidth,
title={Error in minimizer},
xlabel={Average number of quadrature points},
ylabel={Error},
ymajorgrids=true,
grid style=dashed,
legend pos=south west
]
\addplot[dashed,blue,mark=*,mark options={solid,fill=blue}]
    coordinates {
(2143,0.00556530394) (6959,0.00292956738) (21855,0.00149741378)
};
\addplot[dashed,red,mark=square*,mark options={solid,fill=red}]
coordinates {    (183,0.146696585) (631,0.0377325453) (2143,0.00274770301) (6959,1.40889181e-05) 

};
\addplot[dashed,brown,mark=pentagon*,mark options={solid,fill=brown}]
    coordinates {
   (4159,0.0489084233) (18943,0.00421004341) (78079,2.39784785e-05)

};
\addplot[dashed,cyan,mark=triangle*,mark options={solid,fill=cyan}]
    coordinates {
    (1,0.0880207104819) (128,0.0391187400509) (2187,0.022004905962) (16384,0.014084942011)
};
\addplot[dashed,magenta,mark=diamond*,mark options={solid,fill=magenta}]
    coordinates {
    (856.5,0.30836342) (2271.85714286,0.01009461) (3785.0,0.00042114) (3957.0,0.000155917201479)
};
\addplot[dashed,black,mark=star]
    coordinates {
    (405.0,0.126717252947) (1029.0,5.64597095052e-05) 
};
\addplot[dashed,orange,mark=heart,mark options={solid,fill=orange}]
    coordinates {
    (1,0.157720374077) (10,0.0439007920702) (100,0.0167862082653) (1000,0.00298008491817) (10000,0.00180611684029)
};
\legend{tra,cc,gp,puretra,dasgcc,dasggp,10mc}

\end{loglogaxis}
\end{tikzpicture}
\caption{\footnotesize {Compute the additive example with $d=50$.}}
\end{figure}

We test a $7$ dimensional problem. We use an asymmetric distribution $\operatorname{beta}(2,3)$ here. We will get the exact solution with only rounding errors if we use a symmetric distribution. This is because the symmetric construction of the sparse grid will lead to the cancellation of the quadrature points pairs. We still use BFGS method as our optimization solver. The computational results are shown in the Figure $\ref{open loop problem}$. We compare the errors of the $9$ different methods. They are product trapezoidal rule, the average of $10$ runs Monte Carlo, three sparse grid method and three dimension adaptive sparse grid. We only record the data when sparse grid method and dimension adaptive sparse grid method start to converge. As can be seen from the figure, the dimension adaptive sparse grid methods converge faster than classical sparse grid methods for all three univariate rules. The results of sparse gird methods are much better than trapezoidal product rule and Monte Carlo method. 

\begin{figure}[htbp]
\centering 
\subfigure[Clenshaw Curtis with $d=6$]{ 
\begin{minipage}[b]{0.40\linewidth} \includegraphics[width=1\linewidth]{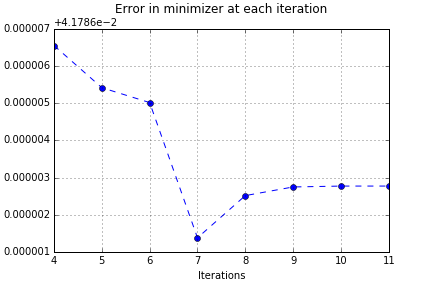}\vspace{4pt} \includegraphics[width=1\linewidth]{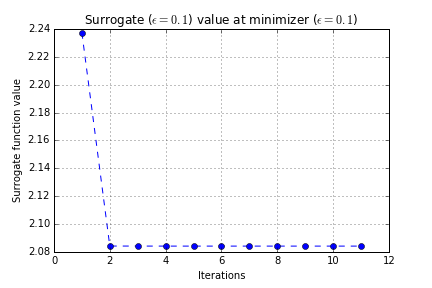}\vspace{4pt} \includegraphics[width=1\linewidth]{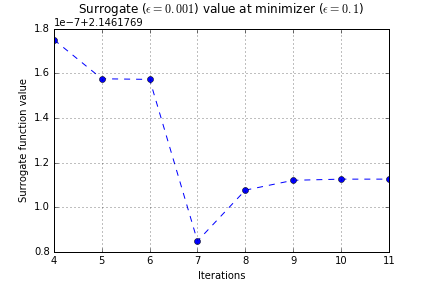}\vspace{4pt} \includegraphics[width=1\linewidth]{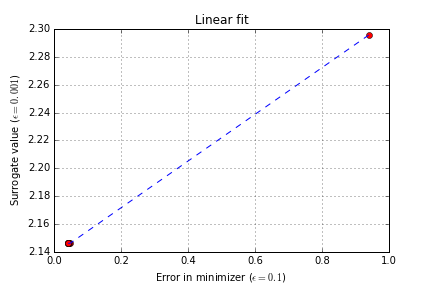} 
\end{minipage}} 
\subfigure[Clenshaw Curtis with $d=7$]{ 
\begin{minipage}[b]{0.40\linewidth} \includegraphics[width=1\linewidth]{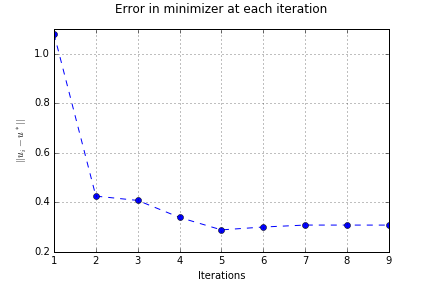}\vspace{4pt} \includegraphics[width=1\linewidth]{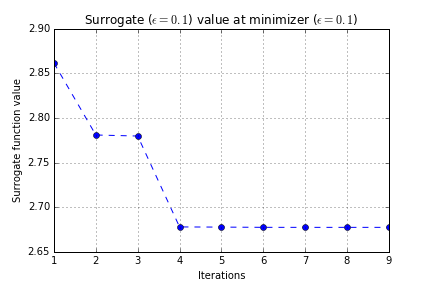}\vspace{4pt} \includegraphics[width=1\linewidth]{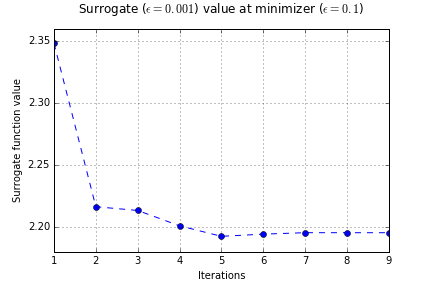}\vspace{4pt} \includegraphics[width=1\linewidth]{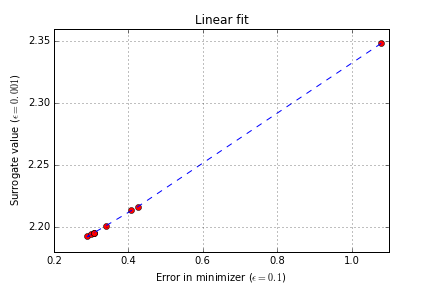} 
\end{minipage}} 

\caption{Solve the stochastic control problem.}\label{solve the sc}
\end{figure}

In Figure $\ref{solve the sc}$, we test our stopping criterion for the stochastic control problem with quadratic cost function and linear dynamic system. For both 6D and 7D examples, our method successfully predicts that we should stop at $7$th iteration for 6D problem and $5$th iteration for 7D problem. 

\section{Conclusions}
We apply the Newton-type methods in solving the stochastic optimisation problem and the dimension-adaptive sparse grid quadrature is used in approximating the integrals involved. In fact we can use more flexible discretisation scheme during the computation if we apply the OTDM. The dimension-adaptive sparse grid quadrature can effectively reduce the computational cost when we use it to compute an integral of which the dimensions are not equally important. When we applied it to solve the stochastic optimisation problem, we find it is more suitable to be used in the OTDM compared with the DTOM. This is because the OTDM allows us to choose 'best' downset in the dimension-adaptive sparse grid formula at each iteration and thus fully exploit the potential of the dimension-adaptive approach. The convergence of the OTDM can be make sure under the condition of Theorem $(\ref{convergence})$. We give the condition when the convergence of our method might break down which leads to a rough stopping criterion. A good stopping criterion is crucial for reducing the computational cost when we solve high dimensional stochastic optimisation problems. We provide another more accurate and practical stopping criterion which only needs reasonable additional computation. We focus on the convex objective function in this paper. For non-convex problems, our approach can only find an approximated local minimizer. In order to solve more general stochastic optimisation problems, other solvers will be taken into consideration in the future research.

\newpage
%\bibliographystyle{plain}

%\bibliography{paper3Bib}
\printbibliography

%\ifx\printbibliography\undefined
%    \bibliographystyle{plain}
%    \bibliography{paper3export}
%\else\printbibliography\fi

\end{document}